\documentclass[preprint,authoryear,12pt]{elsarticle}

\usepackage{graphicx,color}

\usepackage{amsmath,amsthm,amsopn,amstext,amscd,amsfonts,amssymb}

\usepackage{longtable}

\journal{Computational Statistics and Data Analysis}

\theoremstyle{plain} 
\newtheorem{thm}{Theorem}[section]
\newtheorem{prop}[thm]{Proposition}
\newtheorem{cor}[thm]{Corollary}
\newtheorem{remark}[thm]{Remark}

\newcommand{\N}{\mathbb{N}}
\newcommand{\R}{\mathbb{R}}

\newcommand{\restr}[1]{\lower0.4ex\hbox{$|$}\lower0.7ex \hbox{$\scriptstyle{#1}$}}

\begin{document}

\begin{frontmatter}



\title{Maximum likelihood estimation and expectation-maximization algorithm for controlled branching processes}


\author{M. Gonz\'{a}lez}
\author{\corref{cor1}C. Minuesa}
\author{I. del Puerto}

\cortext[cor1]{Corresponding Author: Phone: +34 924289300 ext.
86820. Fax: +34 924272911.\\ E-mail address: cminuesaa@unex.es}
\address{Department of Mathematics, University of Extremadura, 06006 -
Badajoz, Spain.}

\begin{abstract}
  The controlled branching process is a generalization of
  the classical Bienaym\'e-Galton-Watson branching process.
  It is a useful model for describing the evolution of
  populations in which the population size at each
  generation needs to be controlled. The maximum likelihood
  estimation of the parameters of interest for this process
  is addressed under various sample schemes.  Firstly,
  assuming that the entire family tree can be observed, the
  corresponding estimators are obtained and their asymptotic
  properties investigated.  Secondly, since in practice it
  is not usual to observe such a sample, the maximum
  likelihood estimation is initially considered using the
  sample given by the total number of individuals and
  progenitors of each generation, and then using the sample
  given by only the generation sizes.
  Expectation-maximization algorithms are developed to
  address these problems as incomplete data estimation
  problems. The accuracy of the procedures is
  illustrated by means of a simulated example.

\end{abstract}


\end{frontmatter}


\section{Introduction}\label{sec:introdud}

Controlled branching processes are a class of discrete-time
stochastic growth population models characterized by the
existence of a random control mechanism to determine in
each generation (non-overlapping generations) how many
progenitors participate in the subsequent reproduction
process. Once the number of progenitors is known, each one
reproduces independently of the others according to the same
probability law, called the offspring distribution, as usual
in the framework of  branching processes.

In general, the notion of branching has had relevance in the
development of theoretical approaches to problems in such
applied fields as the growth and extinction of populations,
biology (gene amplification, clonal resistance theory of
cancer cells, polymerase chain reactions, etc.),
epidemiology (the evolution of infectious diseases), cell
proliferation kinetics (stem cells, etc.), genetics
(sex-linked genes, mitochondrial DNA, etc.) and algorithm
and data structures (see, for example, the monographs
\cite{Kimmel} and \cite{Haccou}). In particular, the novelty
of adding to the branching notion a mechanism that fixes the
number of progenitors in each generation can allow a great
variety of random migratory movements to be modeled.  The
control mechanism can be defined either by a degenerate
distribution giving rise to deterministic control or in a
random way (through control probability distributions), in
both cases with dependence on the number of individuals in
each generation.  For example, a practical situation that
can be modeled by this kind of process is the evolution of
an animal population that is threatened by the existence of
predators.  In each generation, the survival of each animal
(and therefore the possibility of giving new births) will be
strongly affected by this factor, making the introduction of
a random mechanism (a binomial control process would be reasonable) necessary to
model the evolution of this kind of population.  One can
also model phenomena concerning the introduction or
re-introduction of animal species to inhabit environments in
which they are in potential danger of disappearance or have
previously become extinct.  This re-population can be
achieved by the controlled introduction of new animals until
the species has become firmly established in that habitat.

The family of controlled branching processes includes as
particular cases various models previously introduced in the
branching process literature, such as branching processes with
immigration (see \cite{Sriram-94}), with immigration at
state zero (see \cite{art-Bruss}), with random migration
(see \cite{yanevyanev}), with bounded emigration (see
\cite{libro-ejem}), with adaptive control (see \cite{Bercu}),
and with continuous state space (see \cite{ra}).

The probability theory of this model has been extensively
studied from the pioneering work of \cite{Yanev-75} until the
recent paper of \cite{art-2012} (see also the references therein).
In the last few years, interest in these processes has mainly
focused on the development of their inference theory in order
to guarantee the applicability of these models. Results in this line
 from a frequentist standpoint may be found in
\citet[2005a]{art-2004b} \nocite{art-2005b} for
deterministic control models, using maximum likelihood
estimation, and in \cite{de} and \cite{Sriram} for models
with random control distributions, using martingale
theory (for a multiplicative control function) and weighted
conditional least squares estimation, respectively.

The objective of this paper is to consider the maximum
likelihood estimation of the parameters of interest for a
controlled branching process with random control
distributions under various sample schemes.  Firstly, we
consider the entire family tree until some fixed generation
can be observed.  The results obtained under the observation
of this sample generalize those in
\citet[2005a]{art-2004b}\nocite{art-2005b}. Secondly, since,
in practice, it is not usual to observe the entire family
tree, we consider the maximum likelihood estimation using
initially the sample given by the total number of
individuals and progenitors of each generation, and then the
sample given by only the generation sizes.  We deal with these
problems as incomplete data estimation problems, and develop
expectation-maximization (EM) algorithms to this end (see
\cite{libroEM}, for details of this methodological approach or for recent applications of this methodology in \cite{bzh} and \cite{wwm}).   EM
algorithms have been successfully used to approximate maximum
likelihood estimators when there are missing or incomplete
data, although there are only a few articles on their use in
the context of branching processes (see \cite{vs},
\cite{em-bisexual}, \cite{nina} and \cite{hs}), and in no case for
models which consider random control mechanisms.

After this Introduction, the paper is organized as follows.
We begin by describing the probability model in Section
\ref{sec:model}, in which we introduce some notation and the
working assumptions for the subsequent study. Section
\ref{sec:MLE-complete} is devoted to the maximum likelihood
estimation based on the complete family tree and to studying
the asymptotic properties of the estimators obtained. In
Section \ref{sec:MLE-incomplete}, we address the problem of
obtaining maximum likelihood estimates under incomplete
sampling schemes, developing the EM algorithms. The accuracy
of these algorithms is illustrated by means of a simulated
example in Section \ref{sec:example} (see the supplementary
material for data sets and a further discussion of some
aspects of the example). Some concluding remarks are
provided in Section \ref{sec:conclusions}. Finally, in order
to allow a more readily comprehensible reading, Appendices A, B,
and C are devoted to giving the proofs of the theoretical
results set out in the paper.

\section{The Probability Model}\label{sec:model}

We shall focus our attention on the class of the controlled
branching process with random control function (CBP).
Mathematically, this process is a discrete-time stochastic
growth population model $\{Z_n\}_ {n\geq 0}$ defined
recursively as follows:
\begin{equation}\label{def:model}
Z_0=N,\quad Z_{n+1}=\sum_{j=1}^{\phi_n(Z_{n})}X_{nj},\quad n=0,1,\ldots,
\end{equation}
where $N$ is a non-negative integer, $\{X_{nj}:\
n=0,1,\ldots;j=1,2,\ldots\}$ and
$\{\phi_n(k):n,k=0,1,\ldots\}$ are two independent families
of non-negative integer valued random variables. Also,
$X_{nj}$, $n=0,1,\ldots$, $j=1,2,\ldots$, are independent
and identically distributed (i.i.d.) random variables, and,
for each $n=0,1,\ldots$, $\{\phi_n(k)\}_{k\geq 0}$, are
independent stochastic processes with equal one-dimensional
probability distributions. The empty sum in
\eqref{def:model} is considered to be 0. Let
$p=\{p_k\}_{k\geq 0}$ denote the common probability
distribution of the random variables $X_{nj}$, i.e.,
$p_k=P[X_{nj}=k]$, $k\geq 0$, and $m$ and $\sigma^2$ its
mean and variance (assumed finite), respectively. We also
denote by $\varepsilon(k)=E[\phi_0(k)]$ and
$\sigma^2(k)=Var[\phi_0(k)]$ the mean and the variance of
the control variables (assumed finite too).

Intuitively, $Z_n$ denotes the number of individuals
(particles) in the $n$-th generation and $X_{nj}$ the number
of offspring of the $j$-th individual in the  $n$-th generation.
The probability law $p$ is called the offspring
distribution, and $m$ and $\sigma^2$ are the offspring mean and
variance, respectively. The variable $\phi_n(Z_n)$ represents
a control on the number of progenitors in each generation,
in such a way that when $\phi_n(Z_n)=k$ then $k$ will be the
number of individuals who will take part in the reproduction
process that will determine $Z_{n+1}$. Thus, if $\phi
(Z_n)<Z_n$ then $Z_n - \phi_n(Z_n)$ individuals are removed
from the population (emigration, presence of predators,
etc.), and therefore do not participate in the future
evolution of the process. If $\phi_n(Z_n)>Z_n$ then
$\phi_n(Z_n)- Z_n$ new individuals of the same type are
added to the population (immigration, re-population, etc.). No
control is applied to the population when $\phi_n (Z_n)
=Z_n$. Obviously, if $\phi_n(k)=k$ for all $k$, one obtains
the standard Bienaym\'e-Galton-Watson process.

It is easy to verify that $\{Z_n\}_{n\geq 0}$ is a Markov
chain with stationary transition probabilities. Moreover,
assuming
\begin{enumerate}
  \item [\emph{(a)}] $p_0>0$ or $P[\phi_n(k)=0]>0$, $k>0$,
  \item [\emph{(b)}] $\phi_n(0)=0$ almost surely ($a.s.$),
\end{enumerate}
then 0 is an absorbing state and the states $k=1,2,\ldots$ are transient. Whence it is verified that $P[Z_n\to 0]+P[Z_n\to \infty]=1$.

\vspace*{0.5cm} Let us fix the main parameters of interest
and the working assumptions for the development of their maximum
likelihood estimation. Consider a CBP with an offspring
distribution $p$, whose mean and variance are $m$ and
$\sigma^2$, respectively.
Given that one has different control laws for different
population sizes, the problem of estimating the control
parameters would seem intractable based on samples with a
finite dimension unless the control process is assumed to
have a structure that is  stable over time.
In this sense, formally we
consider CBPs given by
\eqref{def:model} with control distributions belonging to
the power series family of distributions, i.e., for each
$k\geq 0$,
\begin{equation}\label{eq:power-serie-family}
    P[\phi_n(k)=j]=a_k(j)\theta^j A_k(\theta)^{-1},\quad j\geq 0;\theta\in\Theta_k,
\end{equation}
with $a_k(j)$ taking known non-negative values, $A_k(\theta)
=\sum_{j=0}^\infty a_k(j)\theta^j$, and $\Theta_k =
\{\theta>0: 0<A_k(\theta)<\infty\}$ being an open subset of
$\R$. We also assume that the sets $\Theta_k$ are
independent of $k$, so that we shall henceforth drop the index
$k$ from $\Theta_k$, the control parameter space. Moreover,
we assume the following regularity condition:
\begin{equation}\label{eq:regularity-condition}
    \prod_{k\in C} A_k(\theta)=A_{\sum_{k\in C}k}(\theta), \quad\text{ for every } C\subseteq\N;\theta\in\Theta.
\end{equation}

\begin{remark}
  The distribution given in \eqref{eq:power-serie-family} is
  an exponential family which includes many important
  discrete distributions (e.g., Poisson, binomial, negative
  binomial, etc.). The condition
  \eqref{eq:regularity-condition} is a technical hypothesis,
  satisfied by a wide set of probability distributions
  belonging to the exponential family. Hence, the control
  distributions in the model depend on a single parameter
  $\theta$, termed the control parameter, and on the size of
  the population, say $k$.
\end{remark}

\vspace*{0.5cm}

It is well known that:
\begin{eqnarray*}
  \varepsilon(k) &=& \varepsilon(k,\theta) = E[\phi_0(k)]=  \theta \frac{d}{d\theta}\log A_k(\theta),\label{eq:cumulant-med-control}\\
  \sigma^2(k) &=& \sigma^2(k,\theta) = Var[\phi_0(k)]=  \theta \frac{d}{d\theta} \varepsilon(k,\theta).\label{eq:cumulant-var-control}
\end{eqnarray*}

Under condition \eqref{eq:regularity-condition}, it can be
deduced that $\varepsilon(k,\theta)=k \mu(\theta)$, $k\geq
0$, $\theta\in\Theta$, where $\mu(\cdot)$ is a continuous
and invertible function.  From
\eqref{eq:regularity-condition},
$A_k(\theta)=A_1(\theta)^k$, $k\geq 1$, so that
\begin{equation*}\label{eq:expr-epsilon}
    \varepsilon(k,\theta)= \frac{\theta \frac{d}{d\theta}A_k(\theta)}{A_1(\theta)^k}=k\frac{\theta \frac{d}{d\theta}A_1(\theta)}{A_1(\theta)}= k\theta \frac{d}{d\theta}\log(A_1(\theta))=k \varepsilon(1,\theta).
\end{equation*}

Therefore, a family of distributions which verifies
\eqref{eq:regularity-condition} can be re-parametrized making
use of the parameter
$\mu=\mu(\theta)=\varepsilon(1,\theta)$. This parameter can
be termed the  migration parameter because of its intuitive
interpretation: if $\mu<1$, the control law allows one to model
processes with expected emigration; if $\mu > 1$, one can
model processes with expected immigration; and if $\mu=1$,  no migration
is expected.  One also notes that, under
assumption \eqref{eq:regularity-condition},
$\sigma^2(k,\theta)=k\theta\mu'(\theta)$, with $\mu'(\cdot)$
denoting the first derivative of $\mu(\cdot)$.

\vspace*{0.5cm}
\begin{remark}\label{rem:examples1} {Three} interesting
  particular cases of distributions which verify
  \eqref{eq:power-serie-family} and
  \eqref{eq:regularity-condition} are the following:
\begin{enumerate}
\item [(i)] For each $k\geq 0$, take $\phi_n(k)$ to follow a
  Poisson distribution of parameter $k\theta$.
  Consequently, $\mu(\theta)=\theta$.  Hence, depending on
  the value of $\theta$, a CBP with this control function
  can model different migratory processes.  It is easy to
  verify that conditions \eqref{eq:power-serie-family} and
  \eqref{eq:regularity-condition} hold by setting
  $a_k(j)=k^j/j!$ and $A_k(\theta)=e^{k\theta}$.
\item [(ii)] For each $k\geq 0$, take $\phi_n(k)$ to follow
  a binomial distribution of parameters $k$ and $q$.  Taking
  $\theta=q(1-q)^{-1}$, $a_k(j)=\left(
                     \begin{array}{c}
                       k \\
                       j \\
                     \end{array}
                   \right)$, and $A_k(\theta)=(1+\theta)^k$, conditions \eqref{eq:power-serie-family} and \eqref{eq:regularity-condition} can be checked straightforwardly, and $\mu(\theta)=\theta(1+\theta)^{-1}=q$. From a practical viewpoint, this could be a reasonable  control mechanism with which to model situations in which, in each generation, each individual can give birth to offspring in the next generation  with probability $q$, and is removed from the population with probability $1-q$, not participating in its subsequent evolution. As $\mu(\theta)<1$, a CBP with this control distribution always models a case of expected emigration, and, for example, could be useful to model the presence of predators in an animal population.
  \item [(iii)] For each $k\geq 0$, take $\phi_n(k)$ to follow a negative binomial distribution of parameters $k$ and $q$. In this case, conditions \eqref{eq:power-serie-family} and \eqref{eq:regularity-condition} can be checked by setting $\displaystyle{\theta=1-q}$, $a_k(j)=\left(\begin{array}{c}
                       j+k-1 \\
                       j \\
                     \end{array}
                   \right)$, and $A_k(\theta)=(1-\theta)^{-k}$. Moreover, $\mu(\theta)=\theta(1-\theta)^{-1}$. As also was the case for the model considered in (i), this process can model either expected immigration or expected emigration.
\end{enumerate}
\end{remark}

Finally, another parameter of great interest for this family
of processes is what is termed the asymptotic mean growth rate. This
is denoted by $\tau_m$, and is defined in general as
$\lim_{k\to\infty} k^{-1}E[Z_{n+1}|Z_n=k]=\lim_{k\to\infty}
k^{-1}m\varepsilon(k)$ (whenever it exists). Under condition
\eqref{eq:regularity-condition}, $\tau_m=m\mu(\theta)$. This
is the threshold parameter that determines the behaviour of
a CBP in relation to its extinction. Following the
classification of CBPs set out in \cite{art-2005a}, we
shall term a CBP as subcritical, critical, or supercritical depending
on whether $\tau_m$ is less than, equal to, or greater than
unity (emulating the
Bienaym\'e--Galton--Watson process  classification).

In summary, we deal with the problem of estimating $p$, $m$,
$\sigma^2$, $\theta$, $\mu(\theta)$, and $\tau_m$ by making use
of the maximum likelihood estimation based on different
samples.

\section{Maximum Likelihood Estimators with Complete Data}\label{sec:MLE-complete}

In this section, we shall consider the maximum likelihood
estimation of the aforementioned parameters of interest
by assuming that one can observe the entire family tree up to
generation $n$ (complete data), i.e., the random variables
$\{X_{li}: \ 1\leq i\leq \phi_l(Z_l);\ 0\leq l \leq n-1\}$,
or at least $\mathcal{Z}_n^*=\{Z_l(k): 0\leq l\leq n-1;
k\geq 0\}$, where
$Z_l(k)=\sum_{i=1}^{\phi_l(Z_l)}I_{\{X_{li}=k\}}$, $0\leq
l\leq n-1$, $k\geq 0$, with $I_A$ standing for the indicator
function of the set $A$. Intuitively, $Z_l(k)$ represents
the number of individuals in generation $l$ who have
exactly $k$ offspring. It is easily deduced that
$\phi_l(Z_l)=\sum_{k=0}^\infty Z_l(k)$ and
$Z_{l+1}=\sum_{k=0}^\infty kZ_l(k)$, $l=0,\ldots ,n-1$.


\vspace*{0.5cm}
Let $Y_l=\sum_{j=0}^l Z_j$, $\Delta_l=\sum_{j=0}^l \phi_j(Z_j)$, and $Y_{l}(k)=\sum_{j=0}^{l} Z_j(k)$, $l\geq 0$, $k\geq 0$. Intuitively, $Y_l$ and
 $\Delta_l$  denote the total number of individuals and the total number of parents until the $l$-th generation, respectively, and $Y_{l}(k)$ represents the accumulated number up to generation $l$ of individuals who have exactly $k$ offspring. The results presented in this section generalize those given in \citet[2005a]{art-2004b} \nocite{art-2005b} for CBPs with a deterministic control function.

\begin{thm}\label{thm:MLE-complete}
Let $\{Z_n\}_{n\geq 0}$ be a CBP verifying \eqref{eq:power-serie-family} and \eqref{eq:regularity-condition}. The maximum likelihood estimators
 (MLEs) of $p_k$, $k\geq 0$, and $\theta$, based on $\mathcal{Z}_n^*$, are, respectively:
$$\widehat{p}_{k,n}=\frac{Y_{n-1}(k)}{\Delta_{n-1}},\  k\geq 0,\quad\mbox{ and }\quad
\widehat{\theta}_n=\mu^{-1}\left(\frac{\Delta_{n-1}}{Y_{n-1}}\right),$$
where $\mu^{-1}(\cdot)$ denotes the inverse of the function $\mu(\cdot)$.
\end{thm}

The proof is given in Appendix A. 

Using this theorem and the invariance of the MLEs under re-parametrization, the following result is immediate:

\begin{cor}\label{cor:MLE-complete-m-var}
  Let $\{Z_n\}_{n\geq 0}$ be a CBP verifying
  \eqref{eq:power-serie-family} and
  \eqref{eq:regularity-condition}. The MLEs of $m$,
  $\sigma^2$, $\mu(\theta)$, and $\tau_m$ based on
  $\mathcal{Z}_n^*$, are, respectively:
\begin{equation*}
  \widehat{m}_n=\frac{Y_n-Z_0}{\Delta_{n-1}},\quad \widehat{\sigma}^2_n=\sum_{k=0}^\infty(k-\widehat{m}_n)^2\widehat{p}_{k,n},\quad \widehat{\mu}_n=\frac{\Delta_{n-1}}{Y_{n-1}}, \mbox{ and} \quad \widehat{\tau}_{m,n}=\frac{Y_n-Z_0}{Y_{n-1}}.
\end{equation*}
\end{cor}

\vspace*{0.5cm}
For simplicity, when the meaning is clear, we shall drop the index $n$ from $\widehat{p}_{k,n}$ and $\widehat{\tau}_{m,n}$ and write simply $\widehat{p}_k$ and $\widehat{\tau}_{m}$.

\begin{remark}\label{rem:MLE-complete}
\begin{enumerate}
\item [(i)] It is worth noting that to obtain the MLE of the
  offspring distribution, $p$, and its associated
  parameters, $m$ and $\sigma^2$, it is not necessary to
  impose the requirement of any knowledge about the control
  distribution. One can thus address this problem in a
  nonparametric framework, obtaining the same estimators for
  these three parameters.
\item [(ii)] The MLEs of $p_k$ and $m$ are intuitively very
  reasonable because we estimate the probability that an
  individual gives rise to $k$ offspring by the relative
  proportion of parents with $k$ offspring, and the offspring
  mean is estimated by the total number of offspring up to
  a certain generation divided by the number of progenitors
  who have generated those offspring.
\item [(iii)] It can be proved that $\widehat{m}_n$,
  $\widehat{\theta}_n$, $\widehat{\mu}_n$, and
  $\widehat{\tau}_{m}$ are also the MLEs of $m$, $\theta$,
  $\mu(\theta)$, and $\tau_m$, respectively, based on the
  sample
  $\{Z_0,...,Z_{n},\phi_0(Z_0),...,\phi_{n-1}(Z_{n-1})\}$
  (see \cite{jager}, Lemma 2.13.2). Moreover,
  $\widehat{\tau}_{m}$ is also the MLE of $\tau_m$ based on
  $\{Z_0,\ldots,Z_n\}$, following similar arguments.
\end{enumerate}
\end{remark}



\subsection[Asymptotic behaviour]{Asymptotic behaviour}\label{subsec:asymp-behaviour}

In order to investigate the asymptotic properties of the
proposed estimators, it will be necessary to make some
working assumptions. To parameters associated
with the offspring distribution, one does not need to assume that
the control variables belong to a power series family of
distributions. Instead, one only needs to assume that the CBP
$\{Z_n\}_{n\geq 0}$ verifies the following conditions:
\begin{equation}\label{cond:prob-explosion}
\begin{split}
(a)&\ \mbox{There exists }\tau=\lim_{k\to\infty} \varepsilon(k)k^{-1}<\infty, \mbox{ and the sequence }\{\sigma^2(k)k^{-1}\}_{k\geq 1}\\
&\ \hspace{2ex} \mbox{ is bounded.}\\
(b)&\ \tau_m=\tau m >1, \mbox{ and } Z_0 \mbox{ large enough such that } P[Z_n\rightarrow\infty]>0.\\
(c)&\ \{Z_n\tau_m^{-n}\}_{n\geq 0}\mbox{ converges $a.s.$ to a finite random variable } W \mbox{ such that}\\
&\ \hspace{2ex} P[W>0]>0.\\
(d)&\ \{W > 0\}=\{Z_n\to\infty\}\  a.s.
\end{split}
\end{equation}
%

\begin{remark}\label{rem:growth-rate}
\begin{enumerate}
  \item [(i)] In \cite{art-2002}, conditions are provided that guarantee (b) in \eqref{cond:prob-explosion}. Also, in \cite{art-W-ig-explosion}, conditions are established under which $\{W > 0\}=\{Z_n\to\infty\}$ a.s. is verified.
  \item [(ii)]  It can be proved (see \cite{art-2002}, Theorem 4) that, under condition \eqref{cond:prob-explosion}, on the set $\{Z_n\to\infty\}$ one has that
$$\lim_{n\to\infty}Z_n^{-1}Z_{n+1}=\tau_m\quad \text{ a.s.}$$
\end{enumerate}
\end{remark}


\vspace*{0.5cm}
We shall now establish a preliminary result that will be used in the study of the estimators' asymptotic properties. The proof is omitted because it is a consequence of Remark \ref{rem:growth-rate}(ii) and the Stolz-Ces\`aro Lemma.

\begin{prop}\label{prop:asymp-behaviour}
  Let $\{Z_n\}_{n\geq 0}$ be a CBP verifying the conditions
  given in \eqref{cond:prob-explosion}.  Then, on the set
  $\{Z_n\to\infty\}$, it is verified that:
\begin{enumerate}
  \item [(i)] $\lim_{n\rightarrow\infty} Z_n^{-1}\phi_n(Z_n)=\tau$ a.s.
  \item [(ii)] $\sum_{n=0}^\infty \phi_n(Z_n)^{-1}< \infty$ a.s.
  \item [(iii)] $\lim_{n\rightarrow\infty} Y_n^{-1}Y_{n+1}=\tau_m$ a.s.
  \item [(iv)] $\lim_{n\rightarrow\infty }Y_n^{-1}\Delta_n=\tau$ a.s.
  \item [(v)] $\lim_{n\rightarrow\infty }\Delta_n^{-1}\phi_n(Z_n)=\tau_m^{-1}(\tau_m - 1)$ a.s.
  \item [(vi)] $\lim_{n\rightarrow\infty }\varepsilon(Z_n)^{-1}\phi_n(Z_n)=1$ a.s.
\end{enumerate}
\end{prop}

\vspace*{0.5cm} In the following result, we study asymptotic
properties of the estimators related to the offspring
distribution, i.e., $\widehat{p}_k$, $k\geq 0$,
$\widehat{m}_n$, and $\widehat{\sigma}_n^2$.  For simplicity,
we shall use the notation ${\cal D}=\{ Z_n\to\infty\}$ and $P_{\cal
  D}[\cdot]=P[\cdot\mid {\cal D}]$. The result holds
whether or not conditions \eqref{eq:power-serie-family} and
\eqref{eq:regularity-condition} on the control are satisfied.

\begin{thm}\label{thm:asymp-prop-estim}
Let $\{Z_n\}_{n\geq 0}$ be a CBP verifying \eqref{cond:prob-explosion}. Then it holds that:
\begin{enumerate}
\item[(i)] $\widehat{p}_k$, $\widehat{m}_n$, and $\widehat{\sigma}_n^2$ are strongly consistent for $p_k$, $m$, and $\sigma^2$, respectively, on $\{Z_n\to\infty\}$.
\item[(ii)] If $P'$ is a probability measure dominated by $P_{\cal D}$, then for any $x\in \mathbb{R}$:
\begin{enumerate}
\item[(a)] $\displaystyle\lim_{n\to\infty}P'[(p_k(1-p_k))^{-1/2}\Delta_{n-1}^{1/2}(\widehat{p}_k-p_k)\leq x]= \Phi (x),$
\item[(b)] $\displaystyle\lim_{n\to\infty} P'[\sigma^{-1}\Delta_{n-1}^{1/2}(\widehat{m}_n-m)\leq x ]=\Phi(x),$
\item[(c)] If $E[X_{01}^4]<\infty$, then $\displaystyle \lim_{n\to\infty}P'[Var[(X_{01}-m)^2]^{-1/2}\Delta_{n-1}^{1/2}(\widehat{\sigma}_n^2-\sigma^2)\leq x] = \Phi(x)$,
\end{enumerate}
\noindent with $\Phi (\cdot)$ denoting the standard normal distribution function.
\end{enumerate}
\end{thm}

The proof is given in Appendix B. 

\begin{remark}
  Using the previous theorem and Lemma 2.3 in
  \cite{Guttorp}, it is immediate to prove that \emph{(ii)}
  also holds for $P[\cdot|Z_n>0]$. Then, taking into account
  Theorem \ref{thm:asymp-prop-estim} and Slutsky's Theorem,
  and assuming $Z_n>0$, one can obtain asymptotic confidence
  intervals for the parameters $p$, $m$, and $\sigma^2$.
  Thus, for example, the asymptotic confidence interval for
  $m$ at the $1-\alpha$ level, $0<\alpha<1$, is given by
      $$\left[\widehat{m}_n-z_\alpha \left(\widehat{\sigma}_n^2 \Delta_{n-1}^{-1}\right)^{1/2},\widehat{m}_n+z_\alpha \left(\widehat{\sigma}_n^2 \Delta_{n-1}^{-1}\right)^{1/2}\right],$$
      with $z_\alpha$ being such that $1-\Phi(z_\alpha)=\alpha/2$.
\end{remark}

\vspace*{0.2cm} Considering now the parameters of the
control law, let us recall that if the latter belongs to the power
series family of distributions then
\eqref{cond:prob-explosion}(a) holds trivially, and
$\tau=\mu(\theta)$. Denoting \emph{equal in distribution} by $\stackrel{d}{=}$, one has the
following result:

\begin{thm}\label{thm:asymp-prop-estim-theta-mu}
  Let $\{Z_n\}_{n\geq 0}$ be a CBP verifying
  \eqref{eq:power-serie-family},
  \eqref{eq:regularity-condition}, and
  \eqref{cond:prob-explosion}. Then it holds that:
\begin{enumerate}
\item[(i)] $\widehat{\theta}_n$, $\widehat{\mu}_n$ and $\widehat{\tau}_{m}$ are strongly consistent for $\theta$, $\mu(\theta)$ and $\tau_m$, respectively, on $\{Z_n\to\infty\}$.
\item[(ii)] If, for each $l\geq 0$ and $z\geq 0$,
  $\phi_l(z)\stackrel{d}{=}\sum_{s=1}^z X_s(l,z)$, with
  $\{X_s(l,z): 1\leq s\leq z; z\geq 0; l\geq 0\}$ being i.i.d.
  random variables with mean $\mu(\theta)$ and variance
  $\theta \mu'(\theta)$ then, for any $x\in \mathbb{R}$,
    \begin{enumerate}
    \item[(a)] $\displaystyle\lim_{n\to\infty} P_{\cal D}\left[(\theta \mu '(\theta))^{-1/2}Y_{n-1}^{1/2}\left(\widehat{\mu}_n-\mu(\theta)\right)\leq x \right]=\Phi(x),$
    \item[(b)] $\displaystyle\lim_{n\to\infty} P_{\cal D}\left[(\sigma^2\mu(\theta)+m^2\theta \mu '(\theta))^{-1/2}Y_{n-1}^{1/2}\left(\widehat{\tau}_{m}-\tau_m\right)\leq x \right]=\Phi(x),$
    \end{enumerate}
\noindent with $\Phi (\cdot)$ denoting the standard normal distribution function.
\end{enumerate}
\end{thm}

The proof is given in Appendix C.

\begin{remark}\label{rem:MLE-tau}
\begin{enumerate}
\item [(i)] It is worthy of note that the condition set out in
  Theorem \ref{thm:asymp-prop-estim-theta-mu}(ii) is satisfied by the control distributions introduced in Remark
  \ref{rem:examples1}.
\item [(ii)] Theorem \ref{thm:asymp-prop-estim-theta-mu}
 (ii) also holds for $P[\cdot|Z_n>0]$. Again,
  assuming $Z_n > 0$, from this theorem and Slutsky's Theorem,
  and replacing the values $m$, $\sigma^2$, $\theta$, and
  $\mu'(\theta)$ by $\widehat{m}_n$, $\widehat{\sigma}_n^2$,
  $\widehat{\theta}_n$, and $\mu'(\widehat{\theta}_n)$,
  respectively, one can obtain asymptotic confidence
  intervals for the parameters $\mu(\theta)$ and $\tau_m$ at
  the $1-\alpha$ level, $0<\alpha<1$: {\footnotesize
    $$\left[\widehat{\mu}_n-z_\alpha
      \left(\widehat{\theta}_n
        \mu'(\widehat{\theta}_n)Y_{n-1}^{-1}\right)^{1/2},\widehat{\mu}_n+z_\alpha
      \left(\widehat{\theta}_n
        \mu'(\widehat{\theta}_n)Y_{n-1}^{-1}\right)^{1/2}\right],$$
$$\left[\widehat{\tau}_m-z_\alpha \left((\widehat{\sigma}_n^2\mu(\widehat{\theta}_n)+\widehat{m}_n^2\widehat{\theta}_n \mu'(\widehat{\theta}_n))Y_{n-1}^{-1}\right)^{1/2},\widehat{\tau}_m +z_\alpha \left((\widehat{\sigma}_n^2\mu(\widehat{\theta}_n)+\widehat{m}_n^2\widehat{\theta}_n \mu'(\widehat{\theta}_n))Y_{n-1}^{-1}\right)^{1/2}\right],$$}
where $z_\alpha$ is such that $1-\Phi(z_\alpha)=\alpha/2$.
\item[(iii)] Notice that $\widehat{\tau}_m$ is also a strongly consistent estimator for $\tau_m$ on $\{Z_n\to\infty\}$ for CBPs only verifying (\ref{cond:prob-explosion}).
\end{enumerate}
\end{remark}

\section{Maximum Likelihood Estimators with Incomplete Data}\label{sec:MLE-incomplete}

In the previous section, we obtained the MLE of the
parameters of interest ($p$, $m$, $\sigma^2$, $\theta$,
$\mu(\theta)$, and $\tau_m$) based on the sample
$\mathcal{Z}_n^*$. However, in practice, it might be
difficult to observe the entire family tree or the variables
in $\mathcal{Z}_n^*$.  More realistic would be to
suppose that only the total number of individuals and of
progenitors of each generation are known, or even only the
generation sizes. Notice that, with these two samples,
$\widehat{\tau}_m$ is the MLE of $\tau_m$ (see Remark
\ref{rem:MLE-complete}(iii)). Hence, we shall focus
attention on the rest of the parameters. We shall address the problem of the
maximum likelihood estimation under the aforecited samples as an incomplete data estimation
procedure, making use of the EM algorithm and considering
$\mathcal{Z}_n^*$ as \emph{hidden} variables.  Starting with
an initial probability distribution, $p^{(0)}$, and an
initial value of the control parameter, $\theta^{(0)}$, we
will construct a sequence $\{(p^{(i)},\theta^{(i)})\}_{i\geq
  0}$ that will converge to the MLE of $(p,\theta)$. This
iterative method consists of two alternating steps which are
iterated until convergence: the E and the M steps. In the
E step, the expectation of the complete
log-likelihood is calculated using the distribution of the
unobserved data. The values of the parameters which maximize
this expectation are calculated in the following M step.

\subsection[Based on the sample $\{Z_0,\ldots,Z_n,\phi_0(Z_0),\ldots,\phi_{n-1}(Z_{n-1})\}$]{Maximum likelihood estimators based on the
sample $\{Z_0,\ \ldots,Z_n,\ \phi_0(Z_0)$, $\ldots,$ $\phi_{n-1}(Z_{n-1})\}$}\label{subsec:MLE-phi-tot}
We shall determine the MLE of the main parameters of the
model assuming that only the set of random
variables $\overline{\mathcal{Z}}_n=\{Z_0,\ \ldots,Z_n,\
\phi_0(Z_0)$, $\ldots,$ $\phi_{n-1}(Z_{n-1})\}$ can be observed.

Notice that, in accordance with Remark \ref{rem:MLE-complete}(iii), the MLEs of $m$, $\theta$, and $\mu(\theta)$ based on the sample $\overline{\mathcal{Z}}_n$ are $\widehat{m}_n$, $\widehat{\theta}_n$, and $\widehat{\mu}_n$, respectively. Hence, we shall focus on finding the MLEs of $p$ and $\sigma^2$ based on this sample, although we present the method in a general way, considering all the parameters.

\subsubsection{The E step}\label{subsubsec:E-step-phi}
We shall present the E step of the EM algorithm in the
$(i+1)$-st iteration. For each $i$, let
$p^{(i)}=\{p_k^{(i)}\}_{k\geq 0}$ and $\theta^{(i)}$ be the
probability distribution and the control parameter,
respectively, obtained in the $i$-th iteration, and
$\mathcal{Z}_n^*|(\overline{\mathcal{Z}}_n,\{p^{(i)},\theta^{(i)}\})$
the probability distribution of the random vector
$\mathcal{Z}_n^*$ given the sample
$\overline{\mathcal{Z}}_n$ and the parameters $p^{(i)}$ and
$\theta^{(i)}$. For simplicity, in the following, we shall use the notation $E_i^*[
\cdot ] =
E_{\mathcal{Z}_n^*|(\overline{\mathcal{Z}}_n,\{p^{(i)},\theta^{(i)}\})}[
\cdot ]$.

In the proof of Theorem \ref{thm:MLE-complete} (see Appendix
A), Equation \eqref{eq:log-lik-tree} gives the
log-likelihood function $\ell(p,\theta \ | \mathcal{Z}_n^*,\
\overline{\mathcal{Z}}_n) =\ell(p,\theta \ |
\mathcal{Z}_n^*)$, which depends on the unobserved variables $Z_l(k)$,
$0\leq l\leq n-1$, $k\geq 0$. The
expectation of the log-likelihood with respect to the
distribution
$\mathcal{Z}_n^*|(\overline{\mathcal{Z}}_n,\{p^{(i)},\theta^{(i)}\})$
is:
{\small\begin{equation}\label{expr:esper-log-verosim-phi}
  E_i^*[\ell( p,\theta \  |  \mathcal{Z}_n^*,\ \overline{\mathcal{Z}}_n)] =  \Delta_{n-1}\log \theta - \log (A_{Y_{n-1}}(\theta))+ \sum_{l=0}^{n-1}\sum_{k=0}^\infty  E_i^*[Z_l(k)] \log p_k +  E_i^*\left[K\right].
\end{equation}}



Thus, to obtain the value of the above expectation, one has
to determine the distribution of $\mathcal{Z}_n^*$ given
$\overline{\mathcal{Z}}_n$ when the parameters of the models
are $p^{(i)}$ and $\theta^{(i)}$.  Since the individuals
reproduce independently, and the control distributions are
independent of the offspring distribution, one has that, for
$z_0$, $z_{l+1}$, $\phi_l^*$, $z_l(k)\in\N\cup\{0\}$, $k\geq 0$, $0\leq l\leq n-1$
satisfying the constraints $z_{l+1}=\sum_{k=0}^\infty k
z_l(k)$ and $\phi_l^*=\sum_{k=0}^\infty z_l(k)$,

{\small
\begin{align}\label{expr:prob-cond-EM-phi}
    P\big[ Z_l(k)&=z_l(k), 0\leq l\leq n-1, k\geq 0 \big| Z_0=z_0,Z_{l+1}=z_{l+1},\phi_l(Z_l)=\phi_l^*, 0\leq l\leq n-1 \big] = \nonumber \\
    &=\frac{P\big[\{Z_0=z_0\}\cap \bigcap_{l=0}^{n-1} \{Z_{l+1}=z_{l+1},\phi_l(Z_l)=\phi_l^*, Z_l(k)=z_l(k), k\geq 0\}\big]}{P\big[\{Z_0=z_0\}\cap \bigcap_{l=0}^n \{Z_{l+1}=z_{l+1},\phi_l(Z_l)=\phi_l^*\}\big]}\nonumber\\
    &=\prod_{l=0}^{n-1} \frac{P\big[Z_{l+1}=z_{l+1},\phi_l(Z_l)=\phi_l^*,Z_l(k)=z_l(k), k\geq 0 | Z_l=z_l\big]}{P\big[Z_{l+1}=z_{l+1},\phi_l(Z_l)=\phi_l^* | Z_l=z_l\big]}\nonumber\\
       &= \prod_{l=0}^{n-1} \frac{P\big[\sum_{k=0}^\infty k
Z_l(k)=z_{l+1},\phi_l(Z_l)=\phi_l^*,Z_l(k)=z_l(k), k\geq 0| Z_l=z_l\big]}{P\big[Z_{l+1}=z_{l+1},\phi_l(Z_l)=\phi_l^*| Z_l=z_l\big]}\nonumber\\
        &=\prod_{l=0}^{n-1}\frac{P\big[\phi_l(Z_l)=\phi_l^*, Z_l(k)=z_l(k), k\geq 0| Z_l=z_l\big]}{P\big[Z_{l+1}=z_{l+1},\phi_l(Z_l)=\phi_l^*| Z_l=z_l\big]}\nonumber\\
                   &= \prod_{l=0}^{n-1}\frac{P\big[\phi_l(z_l)=\phi_l^*, \sum_{i=1}^{\phi_l(z_l)}I_{\{X_{li}=k\}}=z_l(k), k\geq 0\big]}{P\big[\sum_{i=1}^{\phi_l^*}X_{li} =z_{l+1},\phi_l(z_l)=\phi_l^*\big]}\nonumber\\
    &= \prod_{l=0}^{n-1}\frac{P\big[\sum_{i=1}^{\phi_l^*}I_{\{X_{li}=k\}}=z_l(k), k\geq 0\big]}{P\big[\sum_{i=1}^{\phi_l^*}X_{li} =z_{l+1}\big]}\nonumber\\
    &= \prod_{l=0}^{n-1}\frac{1}{P\big[\sum_{i=1}^{\phi_l^*}X_{li} =z_{l+1} \big]}\cdot \frac{\phi_l^*!}{\prod_{k=0}^\infty z_l(k)! }  \prod_{k=0}^\infty p_k^{(i)z_l(k)}.
\end{align}}

Notice that, although the cardinality of the support of the
reproduction law may be infinite, for each $0\leq l\leq n-1$, once
$z_{l+1}$ and $\phi_l^*$ are known, since
$z_{l+1}=\sum_{k=0}^\infty k z_l(k)$ and
$\phi_l^*=\sum_{k=0}^\infty z_l(k)$, only a finite number of
coordinates of the sequence $\{z_l(k): k\geq 0\}$ are non-null.
From (\ref{expr:prob-cond-EM-phi}), it is clear that to obtain the
distribution
$\mathcal{Z}_n^*|(\overline{\mathcal{Z}}_n,\{p^{(i)},\theta^{(i)}\})$,
first it is enough to know the distributions $(Z_l(k), k\geq
0)|(Z_l,\phi_l(Z_l), Z_{l+1}, \{p^{(i)},\theta^{(i)}\})$, for each
$l=0,\ldots, n-1$. Now, given a fixed generation, say $l$,
assuming that $Z_l=z_l$, $Z_{l+1}=z_{l+1}$ and
$\phi_l(z_l)=\phi_l^*$,  it is needed to determine the sample
space of the vector $(Z_l(k), k \geq 0)$ taking into account that
its possible values $(z_l(k),\ k\geq 0)$ must verify the
constrains $z_{l+1}=\sum_{k=0}^\infty k z_l(k)$ and
$\phi_l^*=\sum_{k=0}^\infty z_l(k)$. After that, their
corresponding probabilities must be obtained following the
equation
$$\frac{1}{P\big[\sum_{i=1}^{\phi_l^*}X_{li} =z_{l+1} \big]}\cdot \frac{\phi_l^*!}{\prod_{k=0}^\infty z_l(k)! }  \prod_{k=0}^\infty p_k^{(i)z_l(k)}.$$
To this end, it is enough to calculate them from a multinomial distribution  of parameters $\phi_l^*$ and $p^{(i)}$ and normalize the obtained probabilities. From this, it is straightforward to obtain the expected values $E_i^*[Z_l(k)]$, $k\geq 0$.
 Notice that this distribution does not
depend on $\theta^{(i)}$ and hence it has no influence on
obtaining $E_i^{*}[Z_l(k)]$.

\subsubsection{The M step}\label{subsubsec:M-step-phi}

In the M step, one calculates the values of the parameters $p$
and $\theta$ which maximize the expectation of the complete
log-likelihood, determined in the previous
step. In other words, one has to find the values
$p^{(i+1)}=\{p_k^{(i+1)}\}_{k\geq 0}$ and $\theta^{(i+1)}$
which maximize the expression
\eqref{expr:esper-log-verosim-phi}, subject to the
constraints $\sum_{k=0}^\infty p_k^{(i+1)}=1$,
$p_k^{(i+1)}\geq 0$, $k\geq 0$.

With a procedure similar to that in the proof of
Theorem \ref{thm:MLE-complete} (see Appendix A) to obtain
the MLEs based on the entire family tree, one obtains
that the values for $p$ and $\theta$ in the $(i+1)$-st
iteration are given by

{\small{$$p_k^{(i+1)} = \frac{\sum_{l=0}^{n-1} E_i^*\left[Z_l(k)\right]}{\sum_{k=0}^\infty \sum_{l=0}^{n-1} E_i^*\left[Z_l(k)\right]}=\frac{\sum_{l=0}^{n-1} E_i^*\left[Z_l(k)\right]}{\sum_{l=0}^{n-1} E_i^*\left[\sum_{k=0}^\infty Z_l(k)\right]}=\frac{\sum_{l=0}^{n-1} E_i^*\left[Z_l(k)\right]}{\Delta_{n-1}},\quad k\geq 0,$$}}
and
$$\theta^{(i+1)}=\mu^{-1}\left(\frac{\Delta_{n-1}}{Y_{n-1}}\right).$$

Intuitively, $p_k^{(i+1)}$ represents the ratio of the
average number (with respect to the probability distribution
determined in the E step) of parents with $k$ offspring to
the total number of progenitors. Notice that
$\theta^{(i+1)}$ does not depend on the iteration $i$
because it is only based on $\overline{\mathcal{Z}}_n$,
which is observed, so that the algorithm reaches the value
$\widehat{\theta}_n$ at the first iteration and then never
leaves it. Hence, as $\theta^{(i)}$ plays no role in
calculating $E_i^*[Z_l(k)]$, at each iteration of the
algorithm based on $\overline{\mathcal{Z}}_n$ only
$p_k^{(i)}$ is updated.  Nonetheless, we include
$\theta^{(i)}$ in the description of the procedure in order
for it to be essentially valid in
both cases considered: when $\overline{\mathcal{Z}}_n$ is
observed and when the sample is only $\{Z_0,\ldots,Z_n\}$
(we shall deal with the latter case in Subsection
\ref{subsec:MLE-tot}).

Indeed, in general, the values $p^{(i+1)}=\{p_k^{(i+1)}\}_{k\geq
0}$ and $\theta^{(i+1)}$ obtained in { the M} step are
used to begin another E step and the process is repeated until the
convergence criterion is satisfied, in which case the process
stops, and the final values are obtained, which we shall denoted
by $\widehat{p}_{n}^{(EM)}=\{\widehat{p}_{k,n}^{(EM)}\}_{k\geq
  0}$ and $\widehat{\theta}_n^{(EM)}$, respectively.
When $\overline{\mathcal{Z}}_n$ is observed,
$\theta^{(i+1)}$ is not needed to begin another $E$ step,
and obviously
$\widehat{\theta}_n^{(EM)}=\widehat{\theta}_n$.

\vspace*{0.5cm} It is straightforward to verify the convergence of the
algorithm by checking the conditions given in \cite{libroEM} on
the continuity and differentiability of the expectation of
the complete log-likelihood function. Consequently, the
sequence $\{(p^{(i)},\theta^{(i)})\}_{i\geq 0}$ converges to
the MLE of $(p,\theta)$ based on the sample
$\overline{\mathcal{Z}}_n$ provided that the likelihood
function $\mathcal{L}(p,\theta| \ \overline{\mathcal{Z}}_n)$
is unimodal.

\vspace*{0.5cm} The EM algorithm also provides the MLE of
$m$, $\sigma^2$, and $\mu(\theta)$ based on
$\overline{\mathcal{Z}}_n$ from the estimates obtained for
the parameter $p$:
$$\widehat{m}_n^{(EM)}=\sum_{k=0}^\infty k \widehat{p}_{k,n}^{(EM)},\quad \widehat{\sigma}_n^{2(EM)}=\sum_{k=0}^\infty \left(k-\widehat{m}_n^{(EM)}\right)^2 \widehat{p}_{k,n}^{(EM)},\quad \widehat{\mu}_n^{(EM)}=\mu(\widehat{\theta}_n^{(EM)}).$$

Obviously, $\widehat{m}_n^{(EM)}=\widehat{m}_n$ and $\widehat{\mu}_n^{(EM)}=\widehat{\mu}_n$. Indeed, for each $i\geq 0$,
$$ m_n^{(i+1)} = \sum_{k=0}^\infty k p_k^{(i+1)} = \frac{\sum_{k=0}^\infty k \sum_{l=0}^{n-1} E_i^*[Z_l(k)]}{\sum_{k=0}^\infty \sum_{l=0}^{n-1} E_i^*[Z_l(k)]} = \frac{Y_n-Z_0}{\Delta_{n-1}}=\widehat{m}_n.$$

\vspace*{0.5cm} In summary, and presented in a general way, the
method to estimate the parameters $p$ and $\theta$, and
consequently $m$, $\sigma^2$, and $\mu(\theta)$, consists
of:
\begin{flushright}
\parbox{13.8cm}{
\begin{enumerate}
  \item [Step 0 ] \emph{$i=0$}. Choose values $\theta^{(0)}$, $0\leq p_k^{(0)}\leq 1$, with $\sum_{k=0}^\infty p_k^{(0)}=1$.
  \item [Step 1 ] \emph{E step}. Based on $p^{(i)}$ and $\theta^{(i)}$
  \begin{enumerate}
    \item [(a)] determine $\mathcal{Z}_n^*|(\overline{\mathcal{Z}}_n,\{p^{(i)},\theta^{(i)}\})$,
    \item [(b)] calculate  $E_i[\ell( p,\theta \  | \mathcal{Z}_n^*,\ \overline{\mathcal{Z}}_n)]$.
  \end{enumerate}
  \item [Step 2 ] \emph{M step}. Calculate the values
  \begin{equation*}
    (p^{(i+1)},\theta^{(i+1)}) = {\arg \max}_{p,\theta} \ E_i[\ell( p,\theta \  | \mathcal{Z}_n^*,\ \overline{\mathcal{Z}}_n)].
  \end{equation*}
    \item [Step 3 ] If $\max \{|p_k^{(i+1)}-p_k^{(i)}|, k\geq 0,\ |\theta^{(i+1)}-\theta^{(i)}|\}$ is less than some convergence criterion, the algorithm halts, and the final values are denoted by $\widehat{p}_n^{(EM)}$ and $\widehat{\theta}_n^{(EM)}$. Otherwise, $i$ is incremented by one unit, and Steps 1-3 are repeated.
  \end{enumerate}}
\end{flushright}

\subsection[Based on the sample $\{Z_0,\ldots,Z_n\}$]{Maximum likelihood estimators based on the sample $\{Z_0,\ldots,Z_n\}$}\label{subsec:MLE-tot}
Now, we shall estimate the parameters with reduced sample
information, assuming that only the total number of
individuals at each generation can be observed. Let us write
$\mathcal{Z}_n=\{Z_0,\ldots,Z_{n}\}$. Although we do not
know exactly what the control function is like or the values
$\phi_0(Z_0),\ldots,\phi_{n-1}(Z_{n-1})$, some information
on the kind of control we are dealing with is necessary, as
will be seen below.

The procedure to obtain the MLE of the model parameters is
almost identical to that of the previous case: making use of
the EM algorithm, one constructs a sequence
$\{p^{(i)},\theta^{(i)}\}_{i\geq 0}$ which will converge to
the MLE of $(p,\theta)$ based on the sample $\mathcal{Z}_n$.

In this case, to determine the expectation of the log-likelihood in the E step, which is
{\small\begin{equation}\label{expr:esper-log-verosim}
  E_i[\ell( p,\theta \  |  \mathcal{Z}_n^*,\ \mathcal{Z}_n)] = E_i[\Delta_{n-1}]\log \theta - \log (A_{Y_{n-1}}(\theta))+ \sum_{l=0}^{n-1}\sum_{k=0}^\infty  E_i[Z_l(k)] \log p_k + E_i[K],
\end{equation}}
\noindent where now $E_i[ \cdot ] =
E_{\mathcal{Z}_n^*|(\mathcal{Z}_n,\{p^{(i)},\theta^{(i)}\})}[
\cdot ]$, one has to know the distribution of
$\mathcal{Z}_n^*$ given $\mathcal{Z}_n$ when the parameters
are $p^{(i)}$ and $\theta^{(i)}$. In this case, it can be
proved that
\begin{align}\label{expr:prob-cond-EM-final}
    P\big[ Z_l(k)&=z_l(k), k\geq 0, l=0,\ldots,n-1 \big| Z_0=z_0,\ldots, Z_{n}=z_{n}\big] = \nonumber\\
    &=\prod_{l=0}^{n-1}\frac{a_{z_l}(\phi_l^*)\theta^{\phi_l^*}A_{z_l}(\theta)^{-1}}{P\big[Z_{l+1}=z_{l+1} | Z_l=z_l\big]}\cdot \frac{\phi_l^*!}{\prod_{k=0}^\infty z_l(k)!}  \prod_{k=0}^\infty p_k^{(i)z_l(k)},
\end{align}
where $z_0$, $z_{l+1}$, $z_l(k)\in\N\cup\{0\}$, $k\geq 0$, $0\leq l\leq n-1$, satisfying
 $\sum_{k=0}^\infty k z_l(k)=z_{l+1}$, and with
$\phi_l^*=\sum_{k=0}^\infty z_l(k)$, $0\leq l\leq n-1$. Equation
(\ref{expr:prob-cond-EM-final}) means that to determine the
distribution
$\mathcal{Z}_n^*|(\mathcal{Z}_n,\{p^{(i)},\theta^{(i)}\})$ is
enough to know the distributions $(Z_l(k),k \geq 0)|(Z_l,
Z_{l+1},\{p^{(i)},\theta^{(i)}\})$, $0\leq l \leq n-1$. Now, for
each fixed $l$, to obtain $(Z_l(k),k \geq 0)|(Z_l,
Z_{l+1},\{p^{(i)},\theta^{(i)}\})$, and given $Z_l=z_l$ and
$Z_{l+1}=z_{l+1}$, first one must consider the possible values for
$\phi_l^*$, determined from the control distribution with
parameters $\theta^{(i)}$ and $z_l$ (notice that, for this
purpose, the kind of control distribution of the process has to be
known). Then, for each fixed value $\phi_l^*$, it is needed to
obtain the sample space of the vector $(Z_l(k), k \geq 0)$ taking
into account that its possible values $(z_l(k), k\geq 0)$ must
verify the constrains $z_{l+1}=\sum_{k=0}^\infty k z_l(k)$ and
$\phi_l^*=\sum_{k=0}^\infty z_l(k)$. Finally their corresponding
probabilities are obtained as the product of probabilities from a
multinomial distribution with parameters $\phi_l^*$ and $p^{(i)}$
by the probability that the control distribution takes the value
$\phi_l^*$ (suitably normalized).

The values of the parameters $p$ and $\theta$ which maximize the expectation of the complete log-likelihood \eqref{expr:esper-log-verosim}, subject to the constraints $\sum_{k=0}^\infty p_k^{(i+1)}=1$, $p_k^{(i+1)}\geq 0$, $k\geq 0$, are:

{\small$$p_k^{(i+1)} = \frac{\sum_{l=0}^{n-1} E_i\left[Z_l(k)\right]}{\sum_{k=0}^\infty \sum_{l=0}^{n-1} E_i\left[Z_l(k)\right]}=\frac{\sum_{l=0}^{n-1} E_i\left[Z_l(k)\right]}{\sum_{l=0}^{n-1} E_i\left[\sum_{k=0}^\infty Z_l(k)\right]}=\frac{\sum_{l=0}^{n-1} E_i\left[Z_l(k)\right]}{E_i\left[\Delta_{n-1}\right]},\quad k\geq 0,$$}
and
$$\theta^{(i+1)}=\mu^{-1}\left(\frac{E_i\left[\Delta_{n-1}\right]}{Y_{n-1}}\right).$$
We shall denote the final values after applying the algorithm to convergence by $\widetilde{p}_{n}^{(EM)}=\{\widetilde{p}_{k,n}^{(EM)}\}_{k\geq 0}$ and $\widetilde{\theta}_n^{(EM)}$, respectively.

Again, it can be checked that the conditions given in
\cite{libroEM} on the continuity and differentiability of
the expectation of the complete log-likelihood function are
satisfied by \eqref{expr:esper-log-verosim}. In this case,
the method also provides estimators for $m$, $\sigma^2$,
and $\mu(\theta)$ based on $\mathcal{Z}_n$:
$$\widetilde{m}_n^{(EM)}=\sum_{k=0}^\infty k \widetilde{p}_{k,n}^{(EM)},\quad \widetilde{\sigma}_n^{2(EM)}=\sum_{k=0}^\infty \left(k-\widetilde{m}_n^{(EM)}\right)^2 \widetilde{p}_{k,n}^{(EM)},\quad \widetilde{\mu}_n^{(EM)}=\mu(\widetilde{\theta}_n^{(EM)}).$$

\section{Simulated Example}\label{sec:example}

We shall illustrate the foregoing results with a simulated example.
We consider a CBP whose offspring distribution is given by
$p_0=0.1084$, $p_1=0.2709$, $p_2=0.3386$, and $p_3=0.2822$,
and the control variables $\phi_n(k)$ follow binomial
distributions with parameters $k$ and $q=0.6$. Thus, the
offspring mean and variance are $m=1.7946$ and
$\sigma^2=0.9443$, respectively; $\theta=1.5$,
$\mu(\theta)=0.6$, and the mean growth rate is
$\tau_m=1.0767$.

In practice, a CBP with  control functions of this kind may
be useful to model the evolution of a population with
the presence of predators. Recall that this binomial control
mechanism models situations in which each individual can
give birth to offspring in the next generation with
probability $q$, or can be removed from the population, not
participating in its subsequent evolution, with probability
$1-q$.

Notice that both $\theta$ and $\mu(\theta)$ determine the
control distribution when the population size is known.
Consequently, we shall focus on the migration
parameter $\mu(\theta)$, which in this case is easily
interpretable.  According to the classification of these
models (see Section \ref{sec:model}), the process considered
in this example is a supercritical CBP with an offspring
mean that is also supercritical, i.e., greater than unity. Notice
that $40\%$ of the individuals do not participate in the
reproduction process for the next generation (i.e., they
are taken by predators).

We simulate the first 30 generations of such a process
starting with $Z_0=1$ individual. The different samples
will be denoted by $z^*_{30}$, $\overline{z}_{30}$, and
$z_{30}$ for that based on the entire family tree, on the individuals and progenitors, and on the population size only, respectively  -- see the supplementary material.
Figure
\ref{im:evol-population} shows the evolution of the number
of individuals and progenitors.  One sees that the
reproduction process makes up for the control process, and,
despite the emigration/predators, the process grows. Thus,
this path seems to belong to the set $\{Z_n\to\infty\}$.
Under the conditions of the example, in \cite{art-2002} and
\cite{art-W-ig-explosion} it is proved that, on the set
$\{Z_n\to\infty\}$, the process grows exponentially with
rate $\tau_m$ (hence, the assumption set out in
\eqref{cond:prob-explosion} holds).

\begin{figure}[h]
 \begin{center}
  \includegraphics[width=5cm]{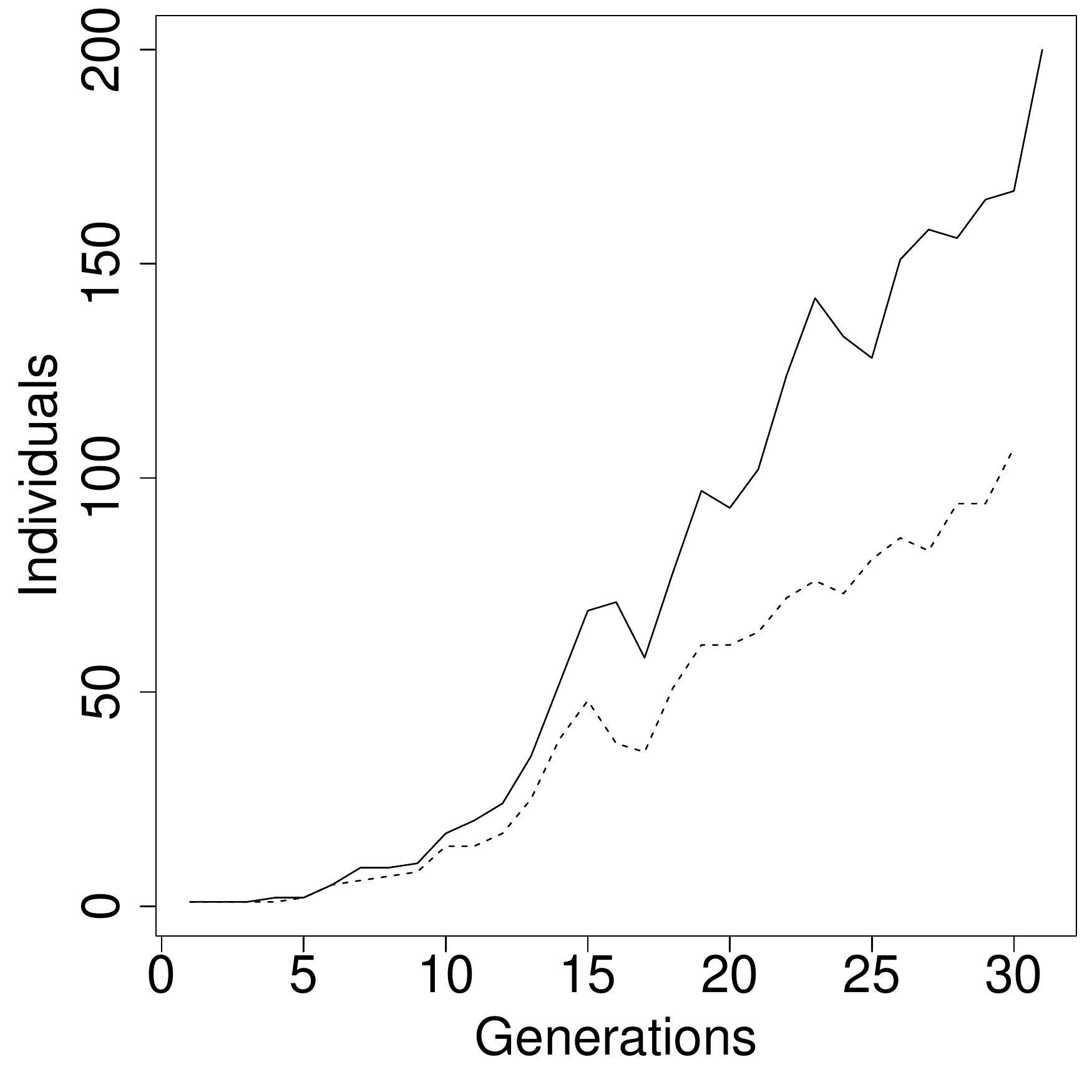}\\
  \caption{\label{im:evol-population}Evolution of the number of individuals (solid line) and the number of progenitors (dashed line).}
 \end{center}
\end{figure}

\begin{figure}[h]
 \begin{center}
  \includegraphics[width=12cm]{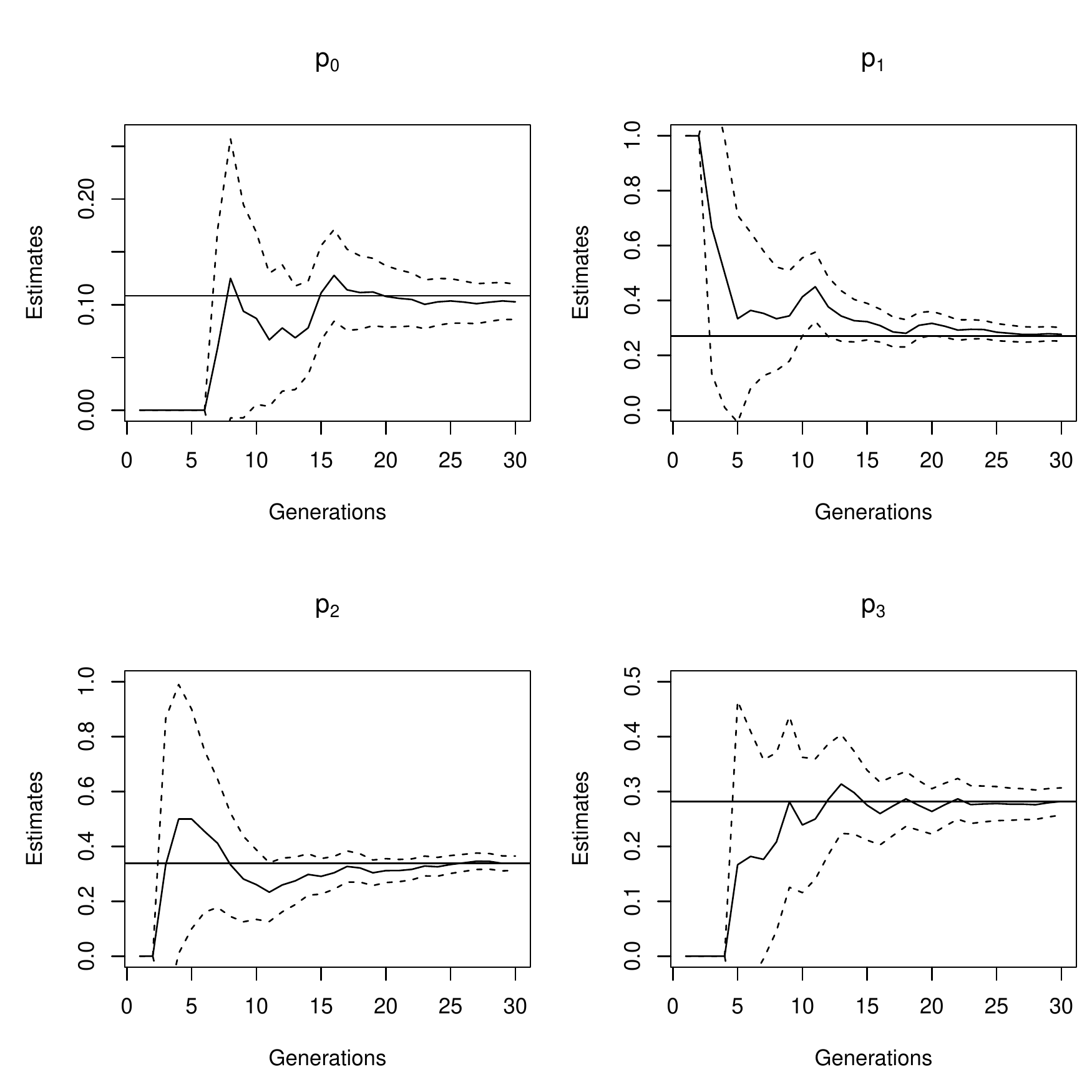}\\
  \caption{\label{im:tree-prob}Evolution of the estimates given by $\widehat{p}_0$, $\widehat{p}_1$, $\widehat{p}_2$, and $\widehat{p}_3$ (solid line), and their approximate $95\%$ confidence intervals (dashed line), together with the true value of the parameters (horizontal line).}
 \end{center}
\end{figure}

\begin{figure}[h]
\begin{center}
\includegraphics*[width=5cm]{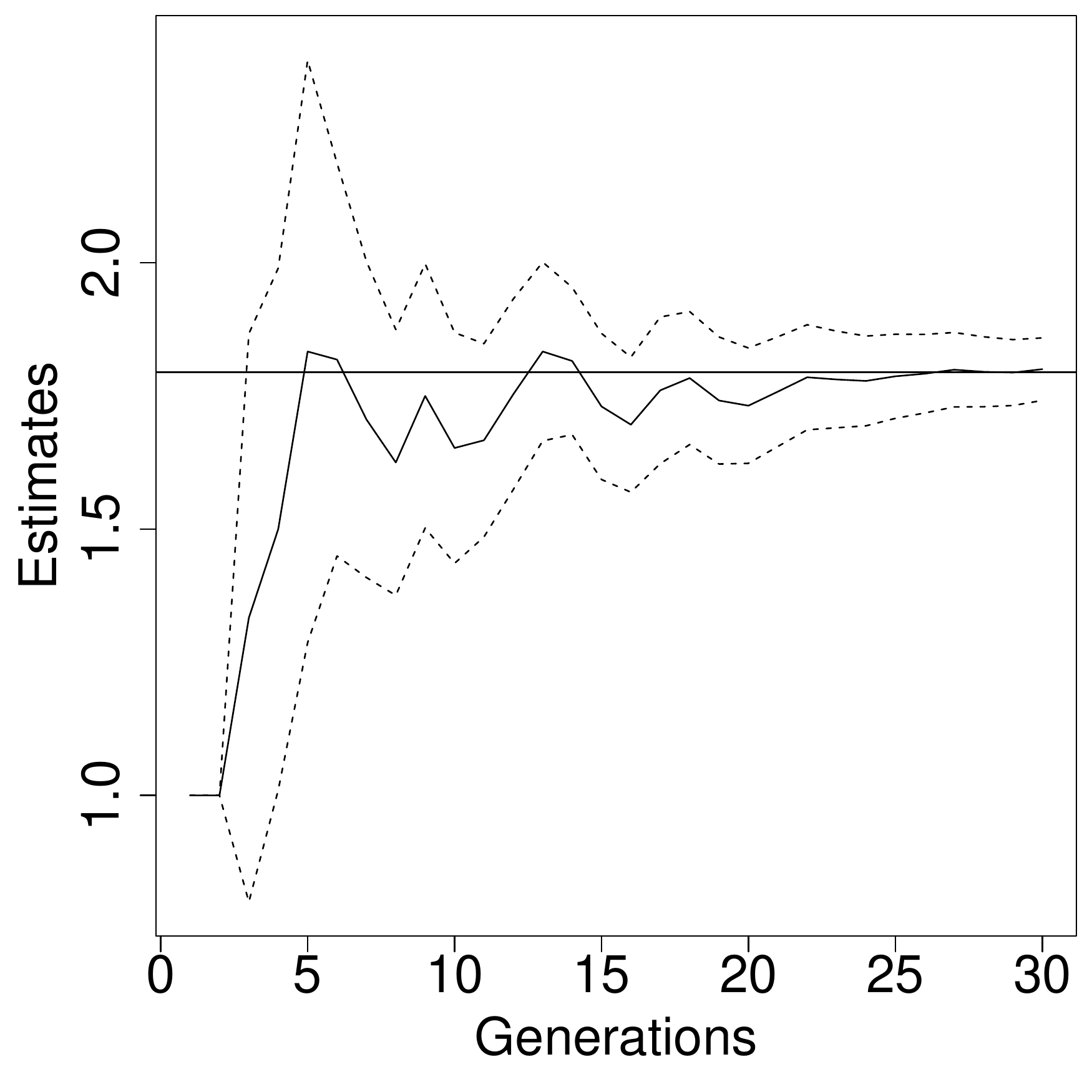}
\hspace{1cm}\includegraphics*[width=5cm]{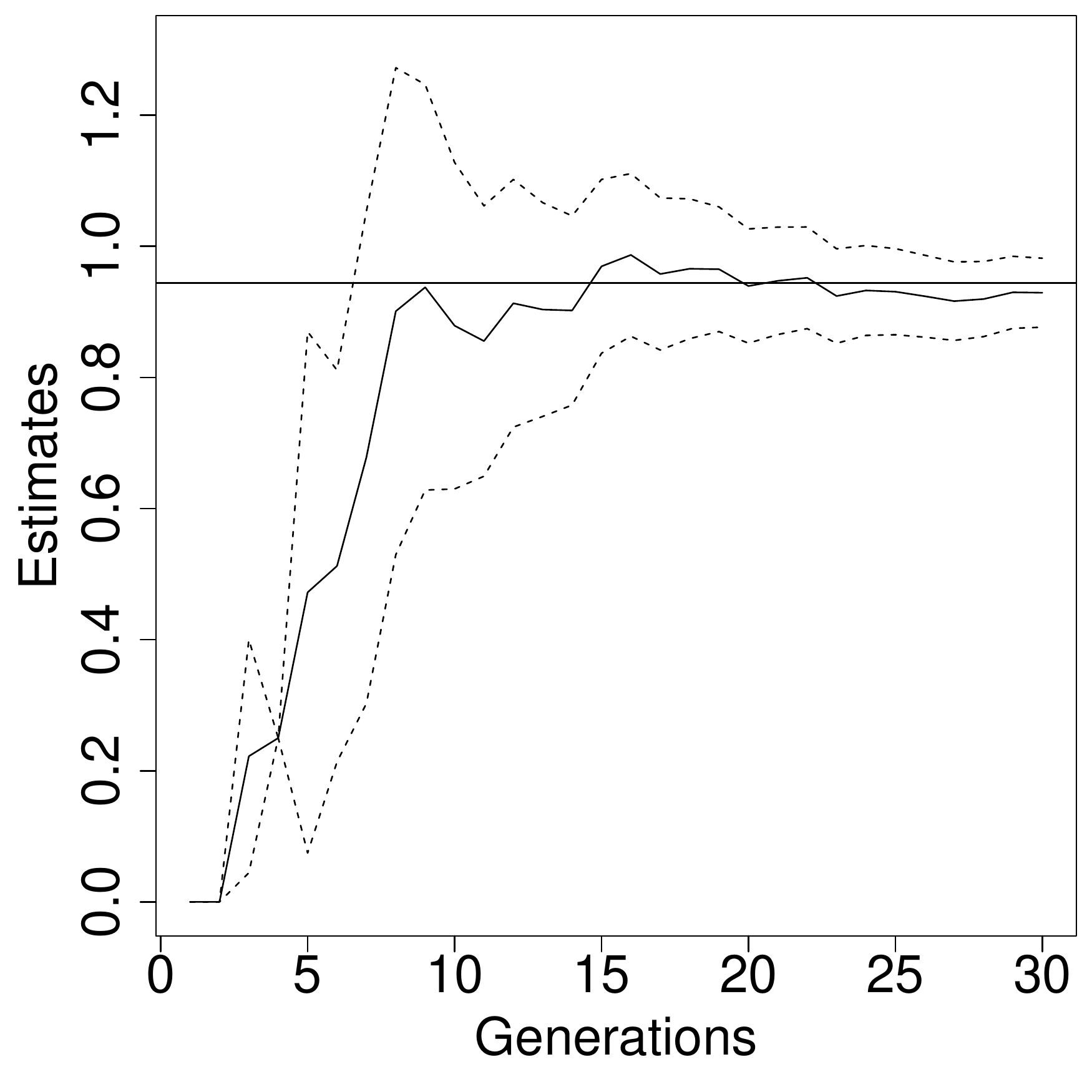}
\caption{\label{im:tree-mean-var} Evolution of the estimates given by $\widehat{m}_n$ (left) and  $\widehat{\sigma}_n^2$ (right) over the course of the generations (solid line) and their approximate $95\%$ confidence intervals (dashed line). The horizontal line represents the true value of each parameter.}
\end{center}
\end{figure}

\begin{figure}[h]
\begin{center}
\includegraphics*[width=5cm]{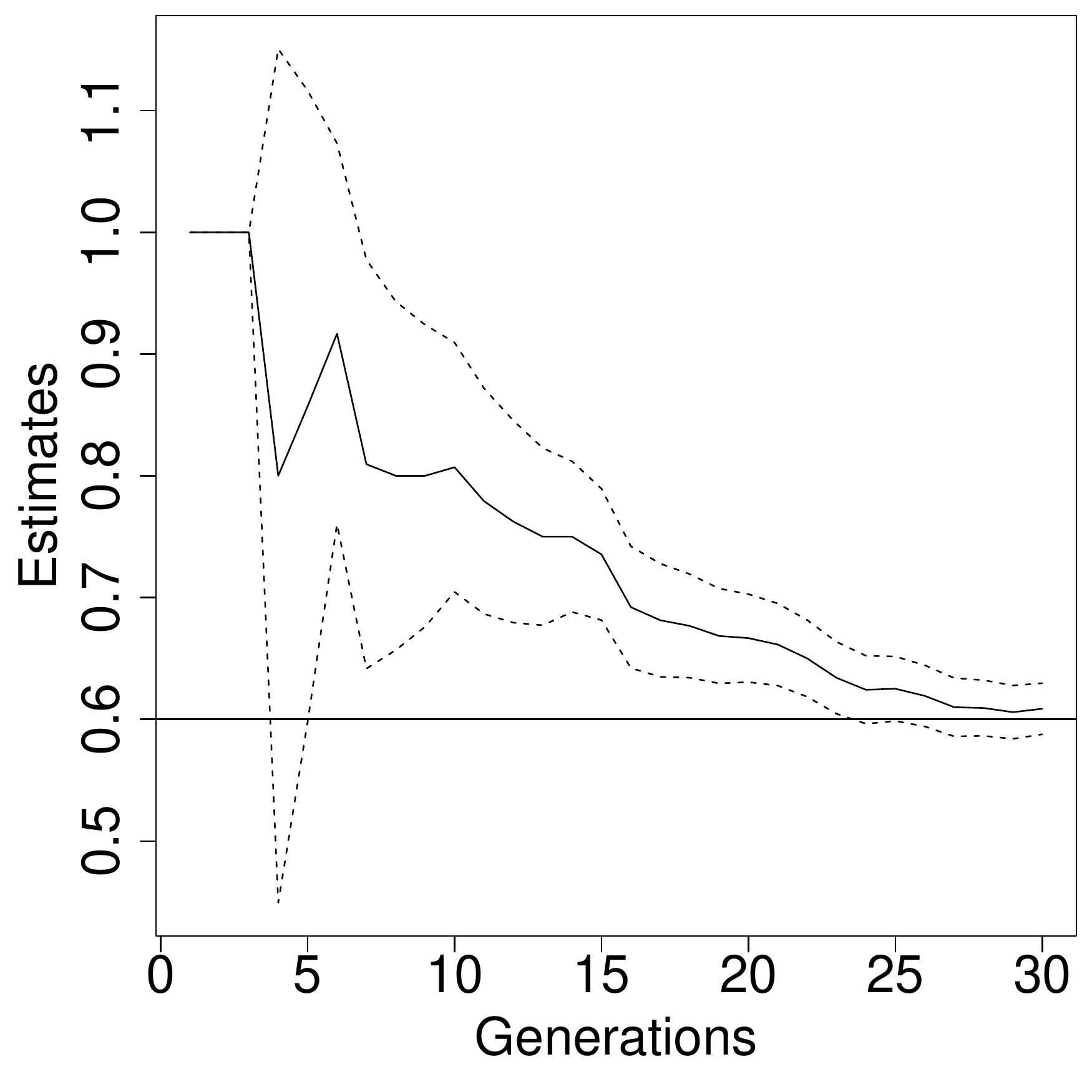}
\hspace{1cm}\includegraphics*[width=5cm]{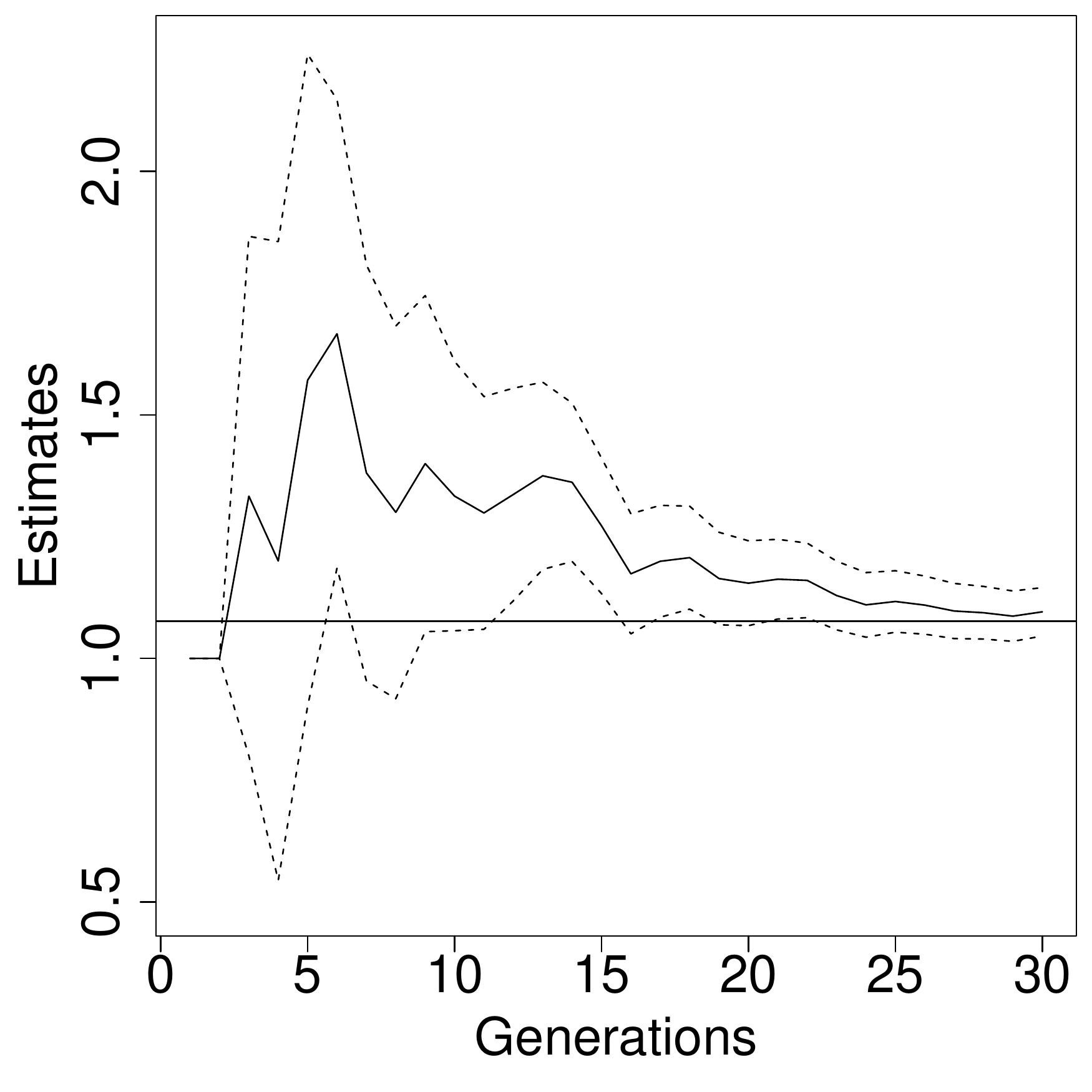}
\caption{\label{im:tree-mu-tau-m} Evolution of the estimates given by $\widehat{\mu}_n$ (left) and $\widehat{\tau}_{m}$ (right) over the course of the generations (solid line) and their approximate $95\%$ confidence intervals (dashed line). The horizontal line represents the true value of each parameter.}
\end{center}
\end{figure}

First, we determined the MLEs and their approximate
$95\%$ confidence intervals based on the entire
family tree, $z^*_{30}$, for $p$, $m$, $\sigma^2$, $\mu(\theta)$, and $\tau_m$. The
estimates are given in Table \ref{t1}.  Figures
\ref{im:tree-prob}--\ref{im:tree-mu-tau-m} show their
behaviours over the course of generations. In these figures we plot the estimates obtained based on the samples restricted to the first $n$ generations, for $n$ going from 0 to 30. One observes that they
approach the true values of the parameters, in accordance with Theorems
\ref{thm:asymp-prop-estim} and \ref{thm:asymp-prop-estim-theta-mu} and Remark
\ref{rem:MLE-tau}(ii).

\begin{table}
  \centering
  \begin{tabular}{ccccccccc}
   \cline{2-9}
\multicolumn{1}{c}{}&\multicolumn{7}{c}{PARAMETERS}\\
\hline
       SAMPLE & $p_0$ & $p_1$ & $p_2$ & $p_3$ & $m$ & $\sigma^2$ & $\mu(\theta)$ & $\tau_m$ \\\hline
    $z_{30}^*$ & .1027 &.2765 & .3389 & .2820 & 1.8002 & .9293 & .6087 & 1.0957 \\
   $\overline{z}_{30}$ & .1211 & .2528 & .3308 & .2953 & 1.8002 & .9927 & .6087 & 1.0957 \\
    $z_{30}$ & .1299 & .3083 & .3283 & .2335 & 1.6653 & .9496 & .6579 & 1.0957 \\\hline
TRUE VALUE     & .1084 & .2709 & .3386 & .2822 & 1.7946 & .9443 & .6000 & 1.0767 \\
      \end{tabular}
  \caption{Estimates of the parameters of interest based on the samples considered $z_{30}^*$, $\overline{z}_{30}$, and $z_{30}$.}\label{t1}
\end{table}


We shall now illustrate the performance of the EM algorithm
in the two situations studied above: using the sample given
by the total number of individuals and progenitors in each
generation, and the sample given by only the generation
sizes. In both cases, assuming that there is no information
available about the offspring distribution, only the maximum
number of offspring per progenitor, we start the algorithm
with a uniform distribution on $\{0,1,2,3\}$ and
$\theta^{(0)}=1/2$.  The maximum number of offspring per
progenitor in an animal population is a datum that is likely
to be known once the reproductive cycle of the corresponding
species is understood.  Even if this information is
unavailable, one can try various (but reasonable) values for this
maximum number of offspring per progenitor, and compare the
results  using the Akaike Information Criterion (AIC) in order
to choose the optimal value (we shall illustrate this procedure
below).

\begin{figure}[h]
 \begin{center}
  \includegraphics[width=12cm]{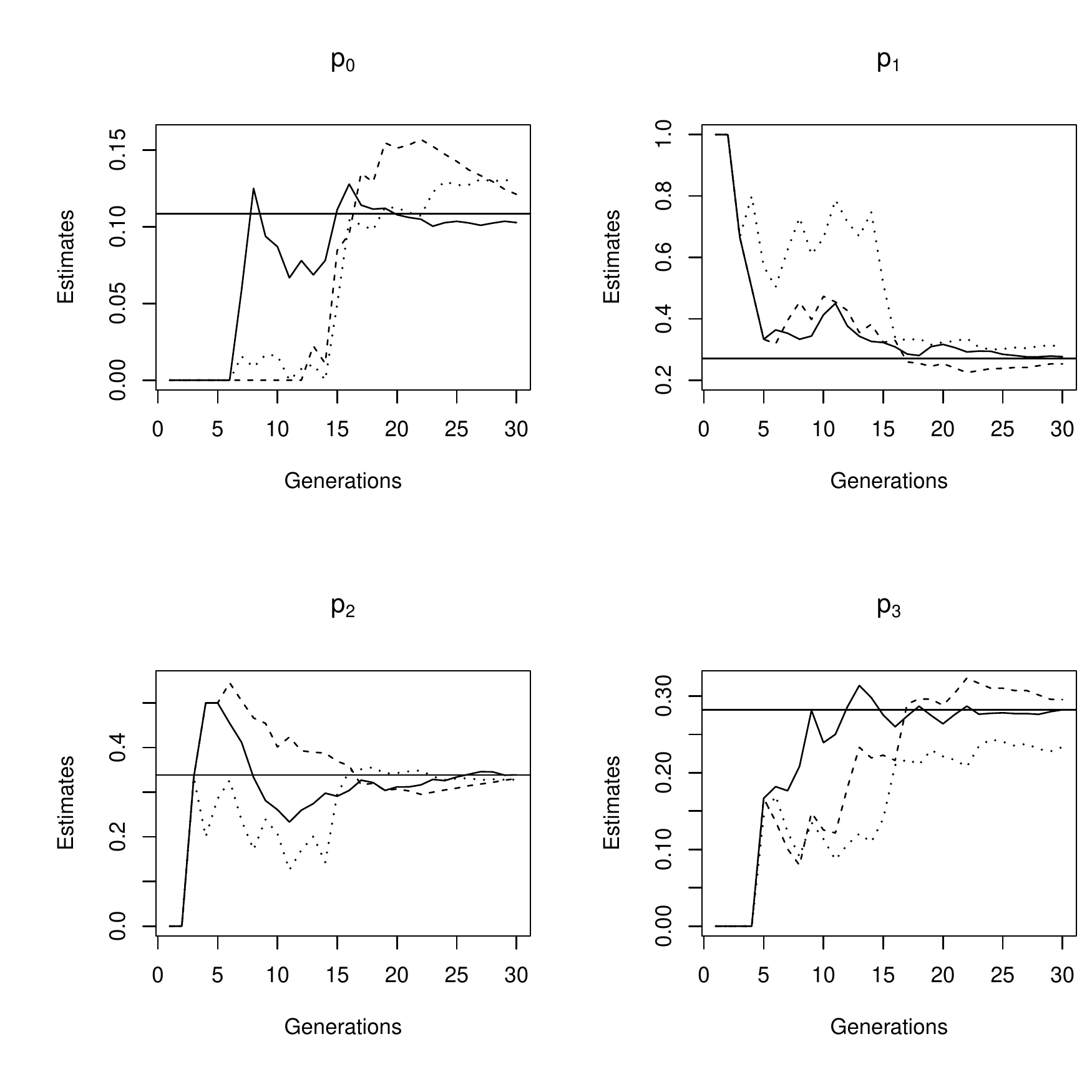}\\
  \caption{\label{im:evol-prob}Evolution of the estimates given by $\widehat{p}$ (solid line), $\widehat{p}^{(EM)}$ (dashed line), and $\widetilde{p}^{(EM)}$ (dotted line).}
 \end{center}
\end{figure}

\begin{figure}[h]
\begin{center}
\includegraphics*[width=5cm]{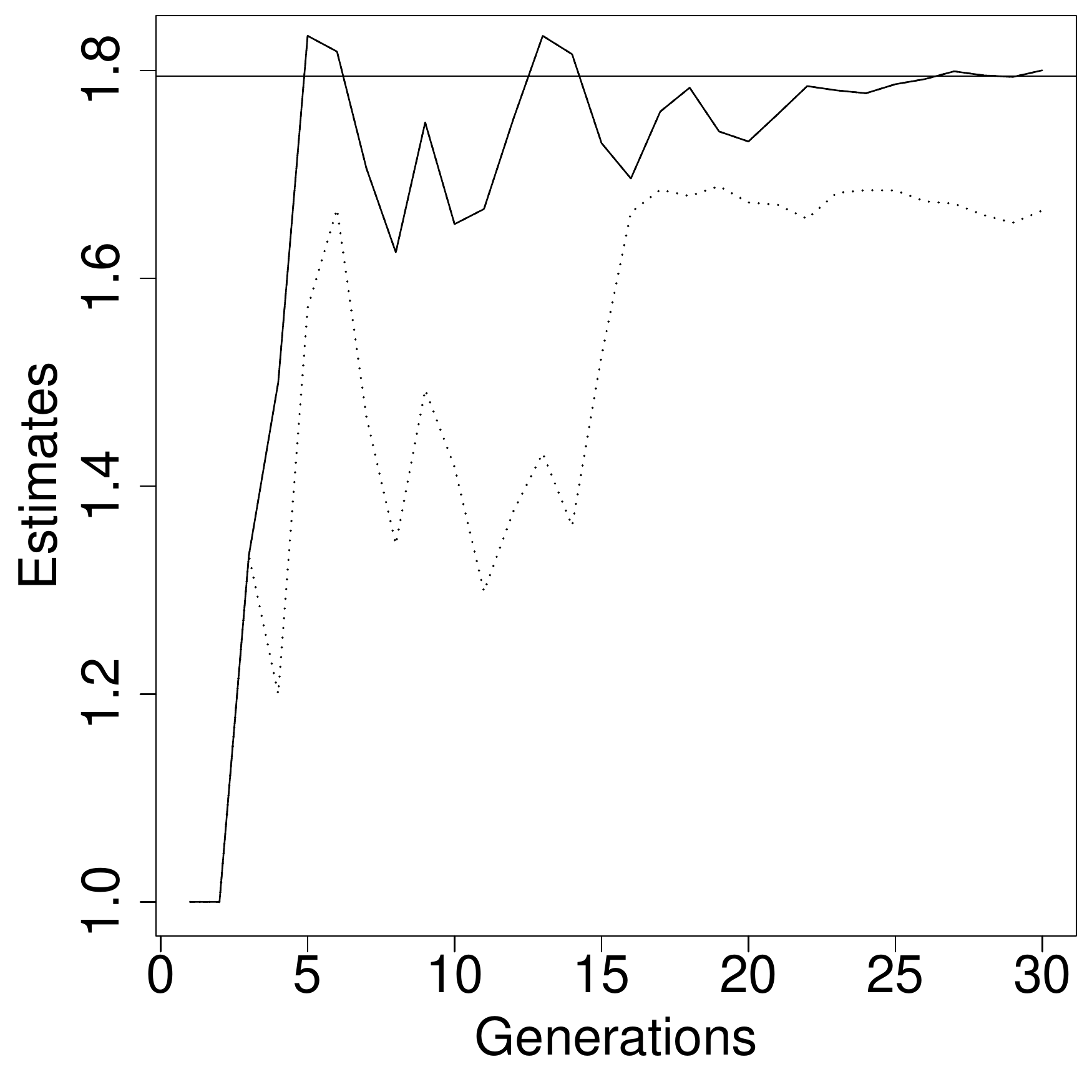}
\hspace{1cm}\includegraphics*[width=5cm]{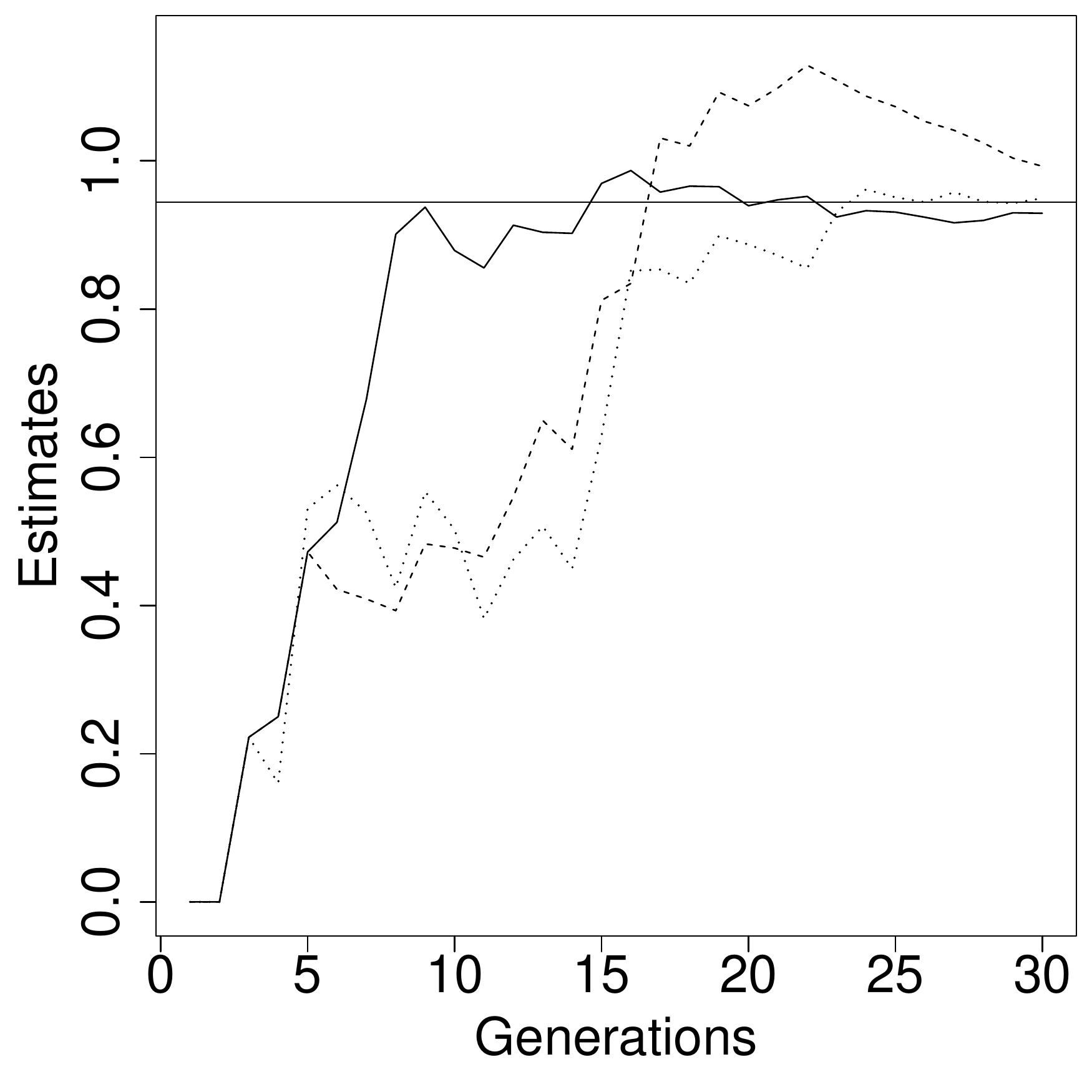}
\caption{\label{im:evol-mean-var} Evolution of the estimates of $m$ (left) and $\sigma^2$ (right) based on the entire family tree (solid line), on the total number of individuals
and progenitors (dashed line) --for estimates of $m$, this coincides with the solid line-- and on the total number of individuals (dotted line), together with the true value of each parameter (horizontal line).}
\end{center}
\end{figure}

Using the first sample, individuals plus progenitors, we ran the
algorithm until attaining a difference between two consecutive
iterations smaller than $10^{-6}$, with this convergence occurring
from iteration 733 onwards. The resulting estimates are given in
Table \ref{t1}. We repeated this procedure considering samples
$\overline{z}_j$, $j=1,\ldots,30$, to assess the consistency of
the estimates. Figures \ref{im:evol-prob} and
\ref{im:evol-mean-var} (right)  show the evolution of these
estimates obtained after convergence of the EM algorithm, and
based on the samples $\overline{z}_j$, $j=1,\ldots,30$ (dashed
lines), together with MLEs based on the entire family tree, for
the parameters $p_k$, $k=0,1,2,3$, and $\sigma^2$.  As was
mentioned above, the EM algorithm is not needed to approximate the
MLEs of $m$, $\theta$, and $\mu(\theta)$ based on the total number
of individuals and progenitors in each generation.

\begin{figure}[h]
\begin{center}
\includegraphics*[width=5cm]{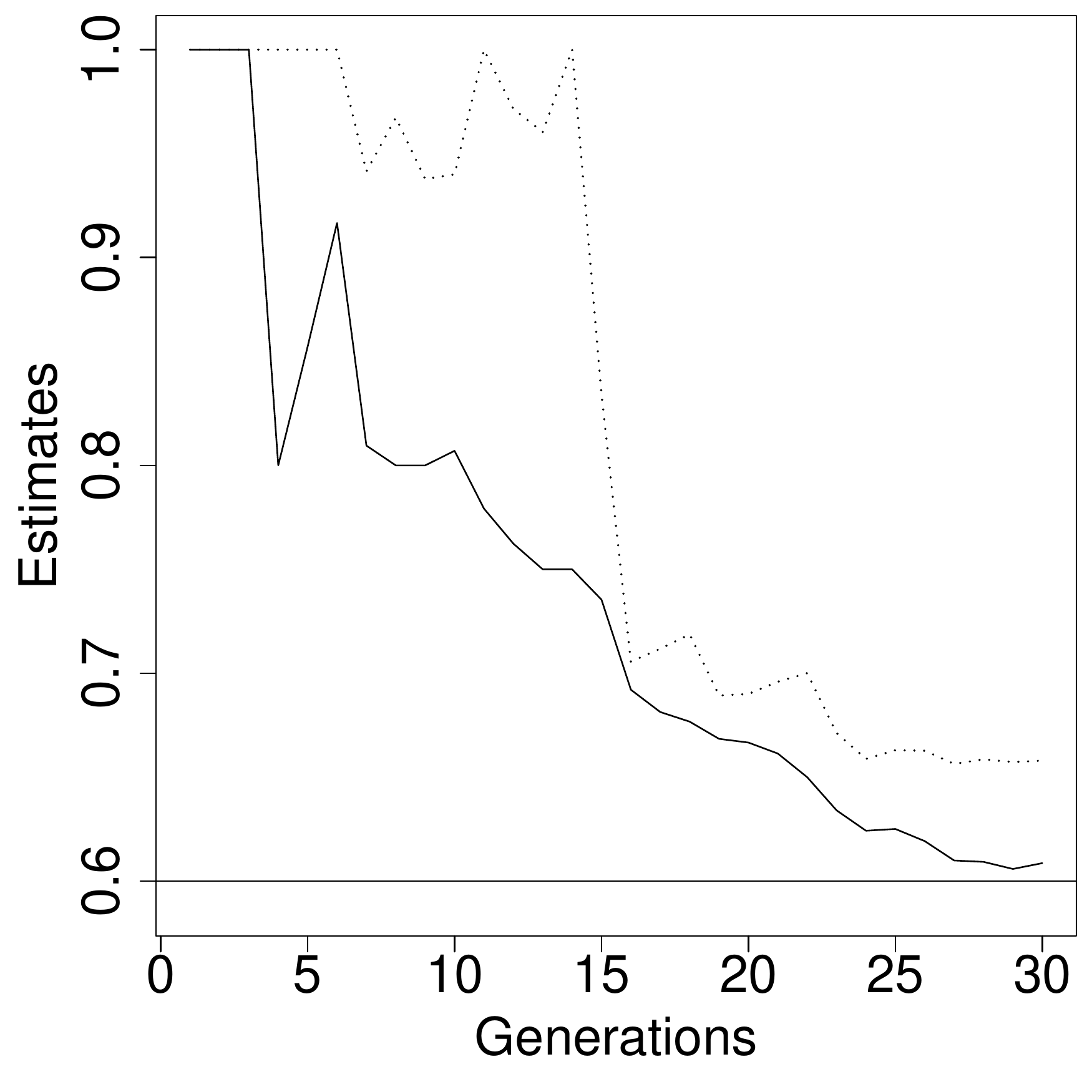}
\hspace{1cm}\includegraphics*[width=5cm]{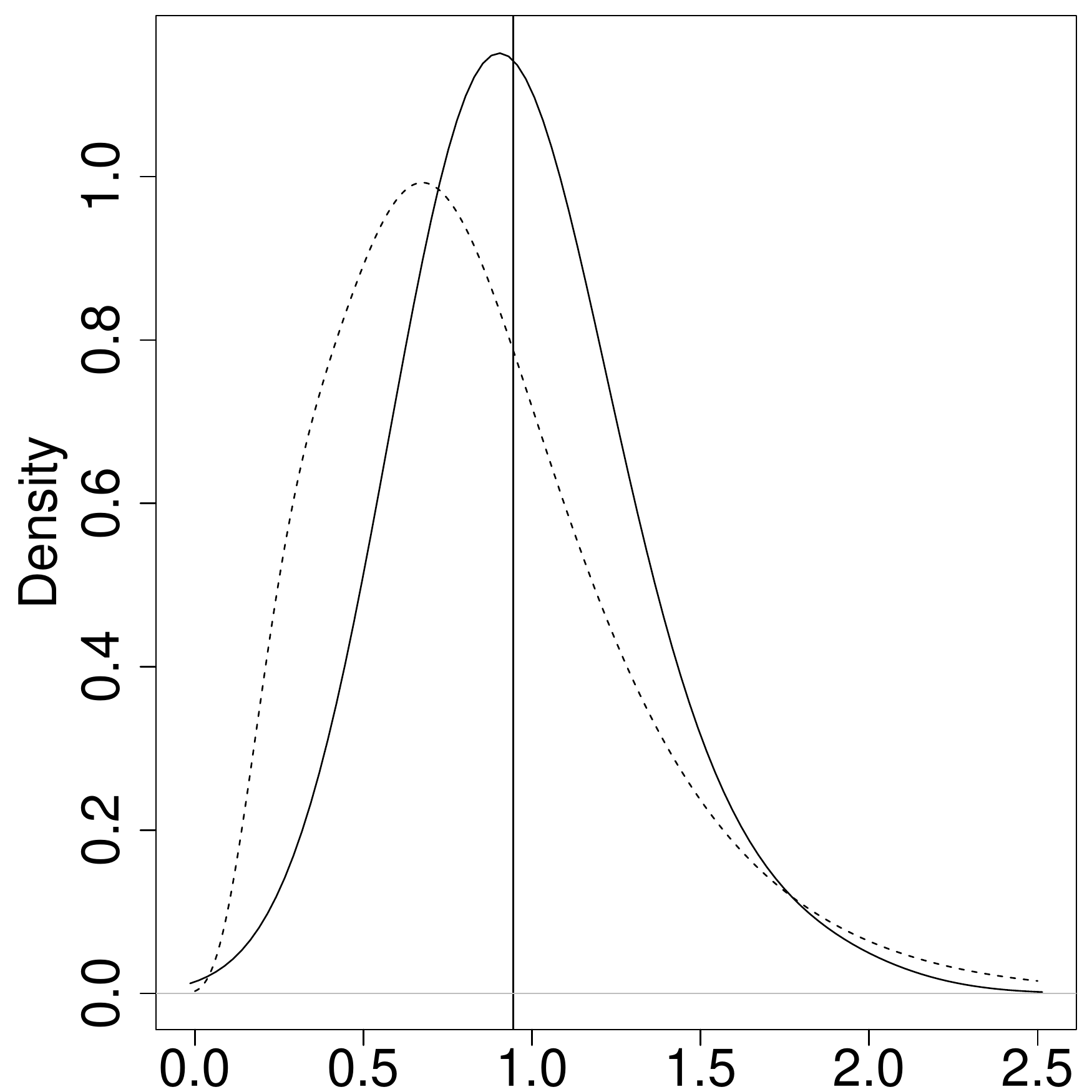}
\caption{\label{im:evol-mu-taum} Evolution of the estimates of $\mu(\theta)$ (left) based on the entire
family tree (solid line) and on the total number of individuals per generation (dotted line), together
 with the true value of the parameter (horizontal line). Bootstrap sampling densities of
$\widehat{\sigma}_{30}^{2(EM)}$ (solid line) and $\widetilde{\sigma}_{30}^{2(EM)}$ (dotted line).}
\end{center}
\end{figure}

We also applied the EM algorithm using the sample defined by
only the population sizes, $z_{30}$. Recall that it is
necessary in this case to know the kind of control
distribution with which one is working. In practice, this
information can come from knowledge of how the population
has developed.  For instance, if there are predators in the
environment, a binomial control distribution would be
clearly justified.  In the simulation, we observed
convergence from iteration 1164 onwards (again for a
precision of $10^{-6}$). The estimates of the parameters are
listed in Table \ref{t1} and their temporal evolution is
plotted in Figures \ref{im:evol-prob},
\ref{im:evol-mean-var}, and \ref{im:evol-mu-taum} (left).
One observes in the figures that all the parameters approach
their respective true values over the course of the
generations.

We studied the influence of the values of $(p^{(0)},\theta^{(0)})$
on the convergence of the algorithms using discrete sensitivity
analysis.  The methods were started with 300 different random
initial values. { The distributions $p^{(0)}$ with
support $\{0,1,2,3\}$ were randomly chosen from a Dirichlet
distribution with all the parameters equal to unity (i.e., by
sampling uniformly from the unit simplex), and the values} of
$\theta^{(0)}$ through the equation
$\theta^{(0)}=q^{(0)}(1-q^{(0)})^{-1}$, with $q^{(0)}$ sampled
from a uniform distribution on the open interval $(0,1)$. Clearly,
the EM algorithm based on the sample $\overline{z}_{30}$ is
insensitive to such choices. But the EM algorithm based on
${z}_{30}$ was observed to not be at all robust to the choice of
initial values, with convergence to different estimates that could
have been local maxima or saddle points. In order to choose the
best approximation to the MLE based on $z_{30}$ (which we will
call the EM estimate), we propose the following methodological
approach.

The log-likelihood function
based on the sample $\mathcal{Z}_n$, denoted by
$\ell(p,\theta\mid \mathcal{Z}_n)$, is given by the
expression
\begin{equation}\label{logexa}
    \ell(p,\theta\mid Z_l=z_l, \ l=0,\ldots, n)=\sum_{j=0}^{n-1}\log\left(\sum_{l=0}^{z_j}P_{z_{j+1}}^{*l}\binom{z_j}{l}
    \frac{\theta^l}{(1+\theta)^{z_j}}\right)
\end{equation}
with $P_{\cdot}^{*l}$ denoting the $l$-fold convolution of the
offspring law {$p$}. While maximization of
$\ell(p,\theta\mid \mathcal{Z}_n)$ would seem to be intractable
using standard methods, (\ref{logexa}) can be evaluated for each
particular $(p,\theta)$. Our proposal is, therefore, to take as EM
estimates of the parameters those associated with the greatest
log-likelihood when it is evaluated at the convergence points of
the EM algorithm started with different randomly chosen values of
the parameters. In our example, the maximum is obtained on the
estimates given in Table \ref{t1} (see the supplementary material
for a further discussion).  {This methodological
strategy can be also followed
when the sample is $\overline{\mathcal{Z}}_n$ (if necessary -not for our  sample observed, $\overline{z}_{30}$), taking into account
that {\small
$$ \ell(p,\theta\mid Z_l=z_l, \phi_l(Z_l)=\phi_l^*, \ l=0,\ldots,
n-1;
Z_n=z_n)=\sum_{j=0}^{n-1}\log\left(P_{z_{j+1}}^{*\phi_l^*}\binom{z_j}{\phi_l^*}
    \frac{\theta^{\phi_l^*}}{(1+\theta)^{z_j}}\right).$$} Moreover,
    it can be extended to processes with whatever type of control
function by only assuming (as has been assumed in our example)
knowledge of the kind of control distribution and of the maximum
number of offspring per progenitor (denoted by $s_{max}$).
Besides, the possibility of calculating the log-likelihood
functions under the samples ${\mathcal{Z}}_n$ and
$\overline{\mathcal{Z}}_n$ allows us to evaluate the influence of
the control distribution and of the value of $s_{max}$, applying
the above method with different control distributions and/or
values of $s_{max}$, and using the AIC to compare the resulting
models.  We have made this study considering the sample
$\overline{z}_{30}$. The results obtained are given in Table
\ref{t3}, in which one observes that for any value of $s_{max}$},
the minimum AIC corresponds to the binomial control distributions.
With respect to the influence of $s_{max}$, the cases $s_{max}=3$
and $4$ led to values that differed little from each other.
Considering therefore \emph{parsimonious parametrization}, it
would be reasonable to choose $s_{max}=3$ as optimal.  In summary,
for problems in which there is no precise knowledge of $s_{max}$
or of the kind of control, a satisfactory procedure would be one
like the foregoing, based on comparing in terms of the AIC several
fitted models (allowing both expected emigration and expected
immigration).

\begin{table}
\begin{center}
{\small
\begin{tabular}{cccccccc}
\cline{1-8}\multicolumn{2}{c}{}&\multicolumn{6}{c}{Control distribution}\\ \cline{3-8}
\multicolumn{1}{c}{$s_{max}$}&\multicolumn{1}{c}{Iterations}
&\multicolumn{2}{c}{Binomial}&
\multicolumn{2}{c}{N Binomial}&\multicolumn{2}{c}{Poisson}\\ \cline{3-4} \cline{5-8}
\multicolumn{1}{c}{}&\multicolumn{1}{c}{}&\multicolumn{1}{c}{$\mbox{Log }$} &
\multicolumn{1}{c}{AIC}&\multicolumn{1}{c}{$\mbox{Log }$} &\multicolumn{1}{c}{AIC}&\multicolumn{1}{c}{$\mbox{Log}$} &
\multicolumn{1}{c}{AIC} \\

\cline{1-8}

3&733& -166.2663& 341.2469& -176.1572& 361.0288& -170.9058& 350.5259\\
4& 4143 &-164.8032 &340.6973 &-174.6942& 360.4792 &-169.4427& 349.9763\\
5& 4244 &-164.8032 &343.1620 &-174.6942 &362.9439& -169.4427 &352.4410\\
6  & 4690  & -164.8032  &  345.7196   &-174.6942 & 365.5015   &-169.4427   &354.9986\\

\cline{1-8}
\end{tabular}}
\caption{Summary of the results for the influence of the control
distributions and $s_{max}$ values. $\mbox{Log }$ denotes
$\ell(\widehat{p}^{(EM)},\widehat{\theta}^{(EM)}\mid
\overline{z}_{30})$. N Binomial denotes the negative binomial distribution. The Iterations column corresponds to the number
of iterations needed to attain a precision of $10^{-6}$ in the
EM procedure for the different $s_{max}$ values.} \label{t3}
\end{center}
\end{table}

\begin{figure}[h]
\begin{center}
\includegraphics*[width=4.25cm]{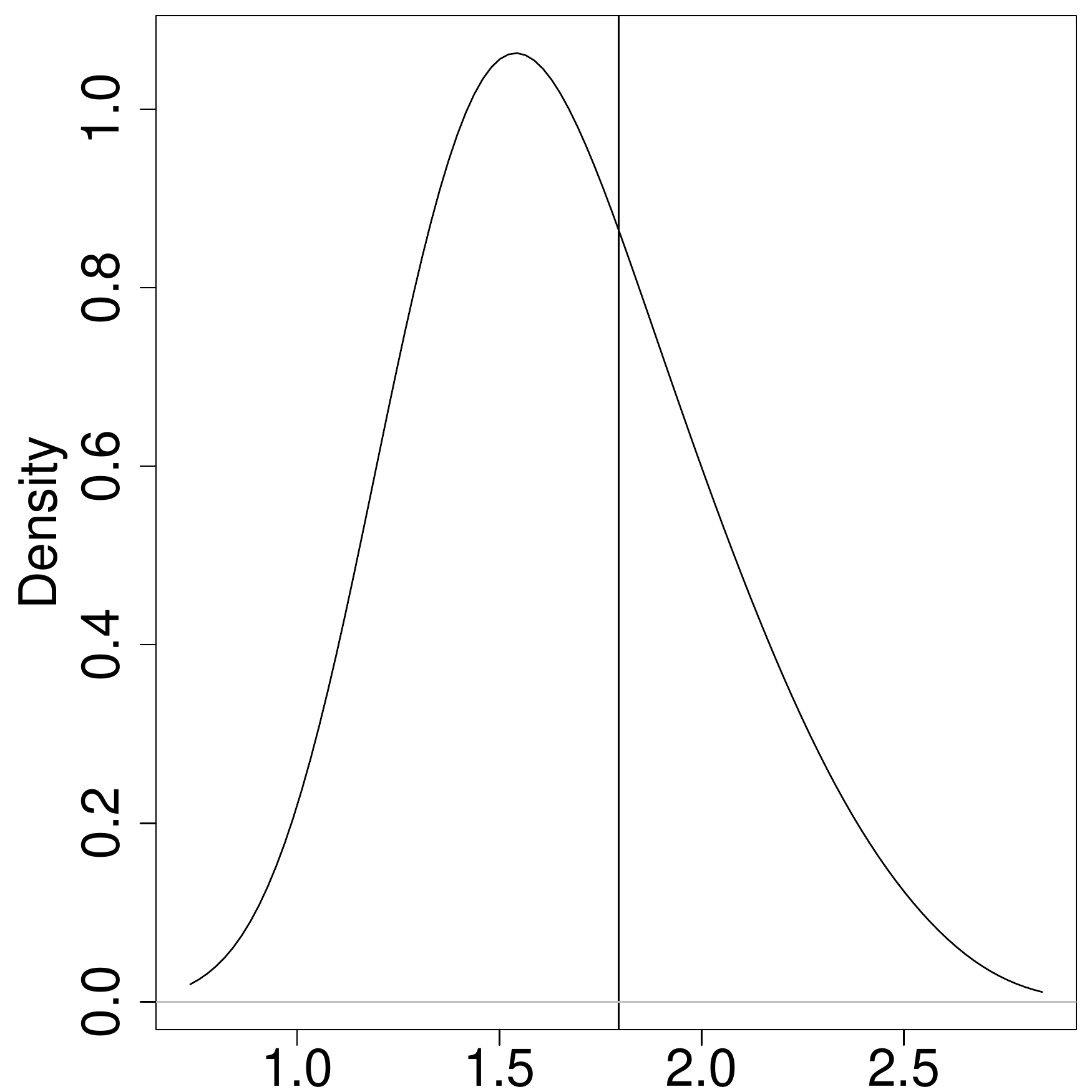}
\hspace{1cm}\includegraphics*[width=4.25cm]{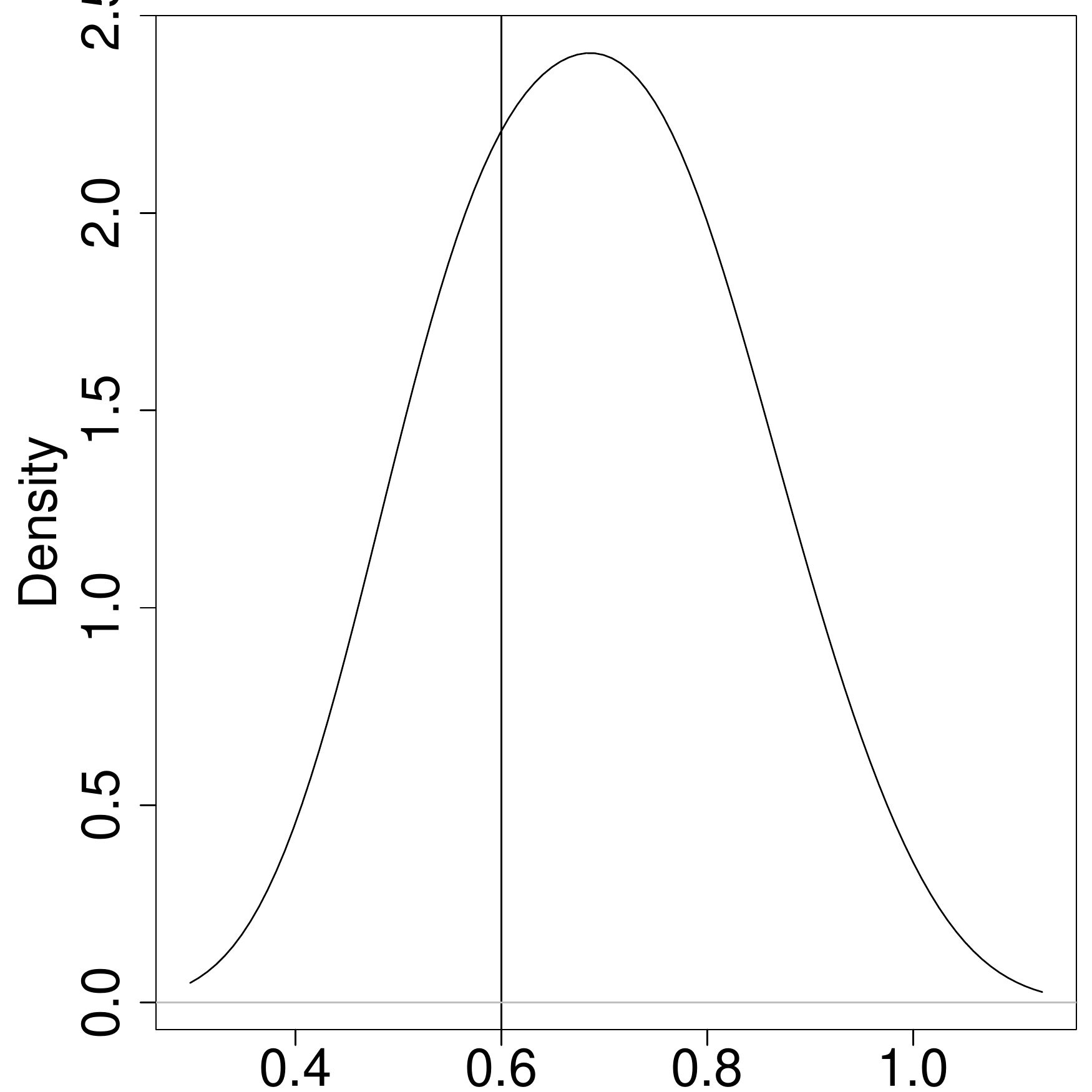}
\hspace{1cm}\includegraphics*[width=4.25cm]{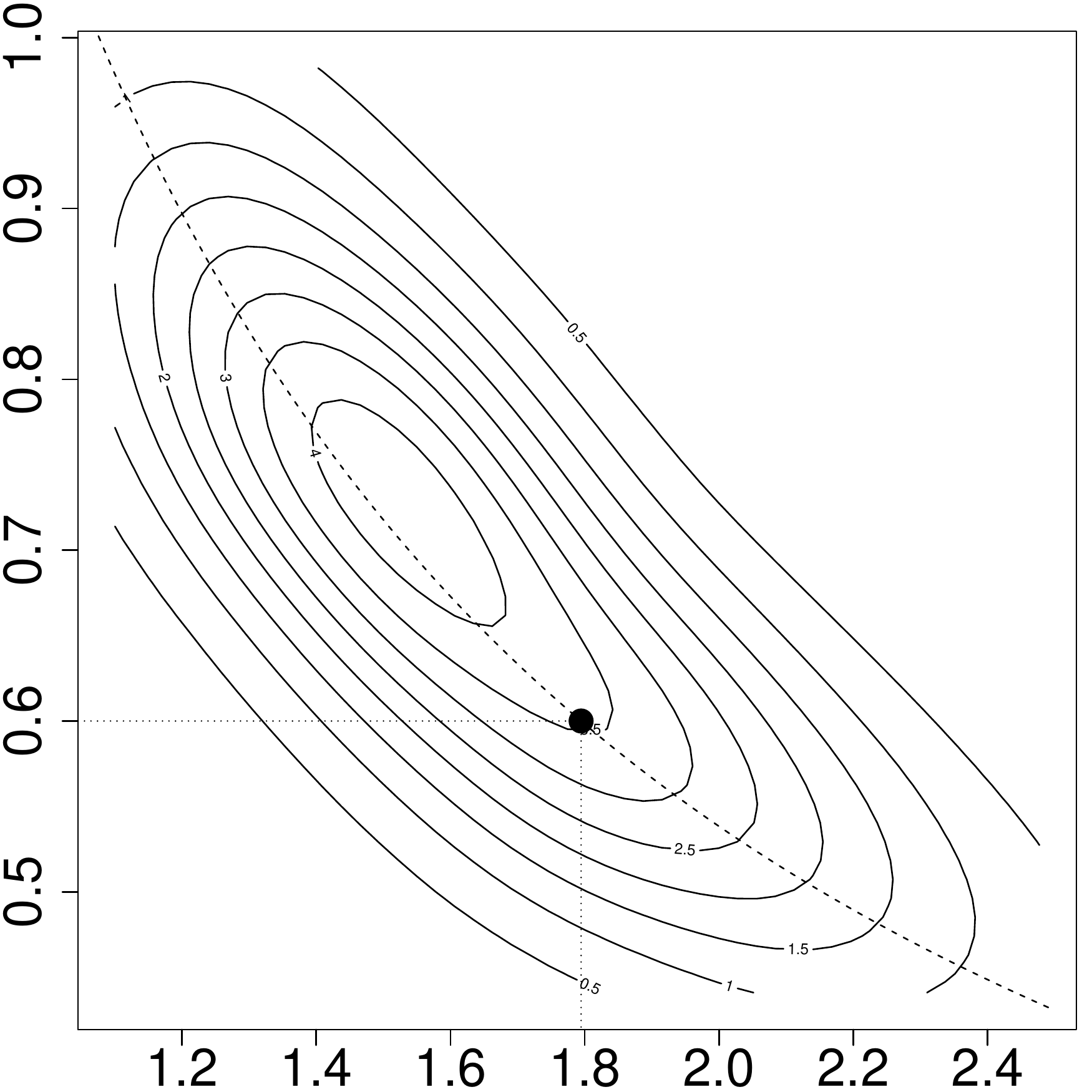}
\caption{\label{new} Bootstrap sampling densities of
$\widetilde{m}_{30}^{(EM)}$ (left) and $\widetilde{\mu}_{30}^{(EM)}$ (center) and joint distribution of
$(\widetilde{m}_{30}^{(EM)},\widetilde{\mu}_{30}^{(EM)})$ with the curve $m\mu(\theta)=1.0767$ (right), together with the true values of the parameters.}
\end{center}
\end{figure}

Finally, to approximate the sampling distributions of
$\widehat{p}_{30}^{(EM)}$, $\widetilde{p}_{30}^{(EM)}$, and
$\widetilde{\theta}_{30}^{(EM)}$ and of their associated
parameters, $\widehat{\sigma}_{30}^{2(EM)}$,
$\widetilde{\sigma}_{30}^{2(EM)}$,
$\widetilde{m}_{30}^{(EM)}$, and
$\widetilde{\mu}_{30}^{(EM)}$, we applied a bootstrap
procedure. We use $\widehat{p}_{30}^{(EM)}$ and
$\widehat{\theta}_{30}^{(EM)}$, based on
$\overline{z}_{30}$, as parameters to perform a Monte Carlo
simulation of 1000 processes up to generation 30.  We
applied the EM algorithm for each of these bootstrapped
samples, obtaining bootstrap approximations to the sampling
distributions of $\widehat{p}_{30}^{(EM)}$, and consequently
of $\widehat{\sigma}_{30}^{2(EM)}$.  Analogously, using the
estimates based on the sample $z_{30}$, we obtained the
bootstrap approximations of the sampling distributions of
the corresponding estimators.  To illustrate these results
without excessive repetition, we shall focus on the
offspring mean and variance and on the migration parameter.
Figure \ref{im:evol-mu-taum} (right) shows the bootstrap
sampling distributions of $\widehat{\sigma}_{30}^{2(EM)}$
and $\widetilde{\sigma}_{30}^{2(EM)}$. One observes that the
distribution of $\widetilde{\sigma}_{30}^{2(EM)}$ is more variable
than that of $\widehat{\sigma}_{30}^{2(EM)}$.  This is a
consequence of the lack of information represented by the control
variables not being observed. Figure \ref{new} shows the
joint distribution of ($\widetilde{m}_{30}^{(EM)}$,
$\widetilde{\mu}_{30}^{(EM)}$) and its marginal
distributions. One observes how strongly these
estimators are related, with their being distributed around the curve given by
$\tau_m=m\mu(\theta)=1.0767$.

Based on the foregoing bootstrap distributions, one can
calculate the bootstrap estimates of the mean squared error
(MSE) of the respective estimators based on the samples
$\overline{z}_{30}$ and $z_{30}$, and compare the accuracy
of the different estimators through their relative
efficiency (eff) (Table \ref{t2}).  One observes from the
table that the estimates provided by assuming observation of
$\bar{z}_{30}$ are preferable to those assuming observation
of $z_{30}$. This can be understood as a reflection of the
greater information content of the former of the two
samples.

\begin{table}
\centering
 \begin{tabular}{cccccccc}
   \hline
   &$p_0$&$p_1$&$p_2$&$p_3$&$m$&$\sigma^2$&$\mu(\theta)$\\\hline
\mbox{MSE based on} $\overline{z}_{30}$  &0.0081&0.0324&0.0332&0.0085&0.0022&0.0808&0.0004\\
\mbox{MSE based on} $z_{30}$ &0.0136&0.0613&0.0637&0.0260&0.1243&0.1438&0.0247\\
\mbox{eff}&1.6779&1.8906&1.9208&3.0594&56.5083&1.7790&62.6719\\
   \hline
 \end{tabular}
\caption{Efficiency of the estimators based on $\overline{z}_{30}$
relative to the estimators based on $z_{30}$ for the parameters of
interest.}\label{t2}
\end{table}

\section*{Computational complexity}

With the aim of determining the order of the computational
complexity of each iteration of the {two  EM algorithms
proposed}, we evaluate the number of operations needed to obtain
$E_i[Z_l(k)]$ and $E^*_i[Z_l(k)]$, $l=0,1,\ldots, n-1$;
{$k=0,1,\ldots, s_{max}$}, respectively (recall
$s_{max}$ is the maximum number of offspring per progenitor).

Let $E^{(i)}=(E_i[Z_l(k)])_{0\leq l\leq n-1; 0\leq k\leq s_{max}}$ and $E^{*{(i)}}=(E^*_i[Z_l(k)])_{0\leq l\leq n-1; 0\leq k\leq s_{max}}$. Considering the sample
$\overline{\mathcal{Z}}_n$, let $B_l$ be the matrix storing the tree associated to the transition from
$\phi_l(Z_{l})$ to $Z_{l+1}$, that is, it stores by rows  the possible
values of the vector $(Z_l(0),\ldots, Z_l(s_{max}))$ such that
$\sum_{k=0}^{s_{max}}Z_l(k)=\phi_l(Z_l)$ and
$\sum_{k=0}^{s_{max}}kZ_l(k)=Z_{l+1}$, $l=0,1,\ldots, n-1$. Let us denote
 $b_l$ the number of rows of $B_l$, $l=0,1,\ldots, n-1$. Finally,  for $l=0,1,\ldots, n-1$, let $P_l$ be a row vector  whose elements are the probabilities of each row of $B_l$, obtained by equation (\ref{expr:prob-cond-EM-phi}), that is, if $\phi_l(Z_l)=\phi^*_l$ and $Z_{l+1}=z_{l+1}$, the corresponding element of $P_l$ for the row of $B_l$ given by $(z_l(0),\ldots, z_l(s_{max}))$ is equal to
 \begin{equation*}
    \frac{1}{P[\sum_{i=1}^{\phi_l^*}x_{li}=z_{l+1}]}\frac{\phi_l^*!}{\prod_{k=0}^{s_{max}}z_l(k)!}\prod_{k=0}^{s_{max}}p_k^{(i)z_l(k)}.
 \end{equation*}
  Then, the $l$-th row of $E^{(i)}$ is equal to the product $P_l\cdot B_l$, $l=0,1,\ldots, n-1$.

  Analogously, assuming the sample $\mathcal{Z}_n$, for $l=0,1,\ldots, n-1$, let $B^*_l$ be the matrix storing the tree associated to the transition from
$Z_{l}$ to $Z_{l+1}$, that is, its rows store  all the possible
values of the vector $(Z_l(0),\ldots, Z_l(s_{max}))$ that allow reaching $Z_{l+1}$ from $Z_l$. To obtain such a matrix, if $Z_l=z_l$ and $Z_{l+1}=z_{l+1}$, we consider for each possible value of $\phi_l(z_l)$, say $\phi_l^*$, every vector $(z_l(0),\ldots, z_l(s_{max}))$ such that $\sum_{k=0}^{s_{max}}z_l(k)=\phi_l^*$ and
$\sum_{k=0}^{s_{max}}kz_l(k)=z_{l+1}$, $l=0,1,\ldots, n-1$. Now, for each one of these vectors we obtain the probabilities (see equation (\ref{expr:prob-cond-EM-final}))

 \begin{equation*}
    \frac{a_{z_l}(\phi_l^*)\theta^{\phi_l^*}A_{z_l}(\theta)^{-1}}{P[Z_{l+1}=z_{l+1}| Z_{l}=z_{l}]}\frac{\phi_l^*!}{\prod_{k=0}^{s_{max}}z_l(k)!}\prod_{k=0}^{s_{max}}p_k^{(i)z_l(k)},
 \end{equation*}
 which are ordered in the row vector $P_l^*$.  Then, the $l$-th row of $E^{*(i)}$ is equal to the product $P^*_l\cdot B^*_l$, $l=0,1,\ldots, n-1$. Let us denote by  $b_l^*$  the number of
the rows of $B_l^*$, $l=0,1,\ldots, n-1$.

\begin{figure}[h]\label{complexity}
 \begin{center}
 \includegraphics[width=5.5cm]{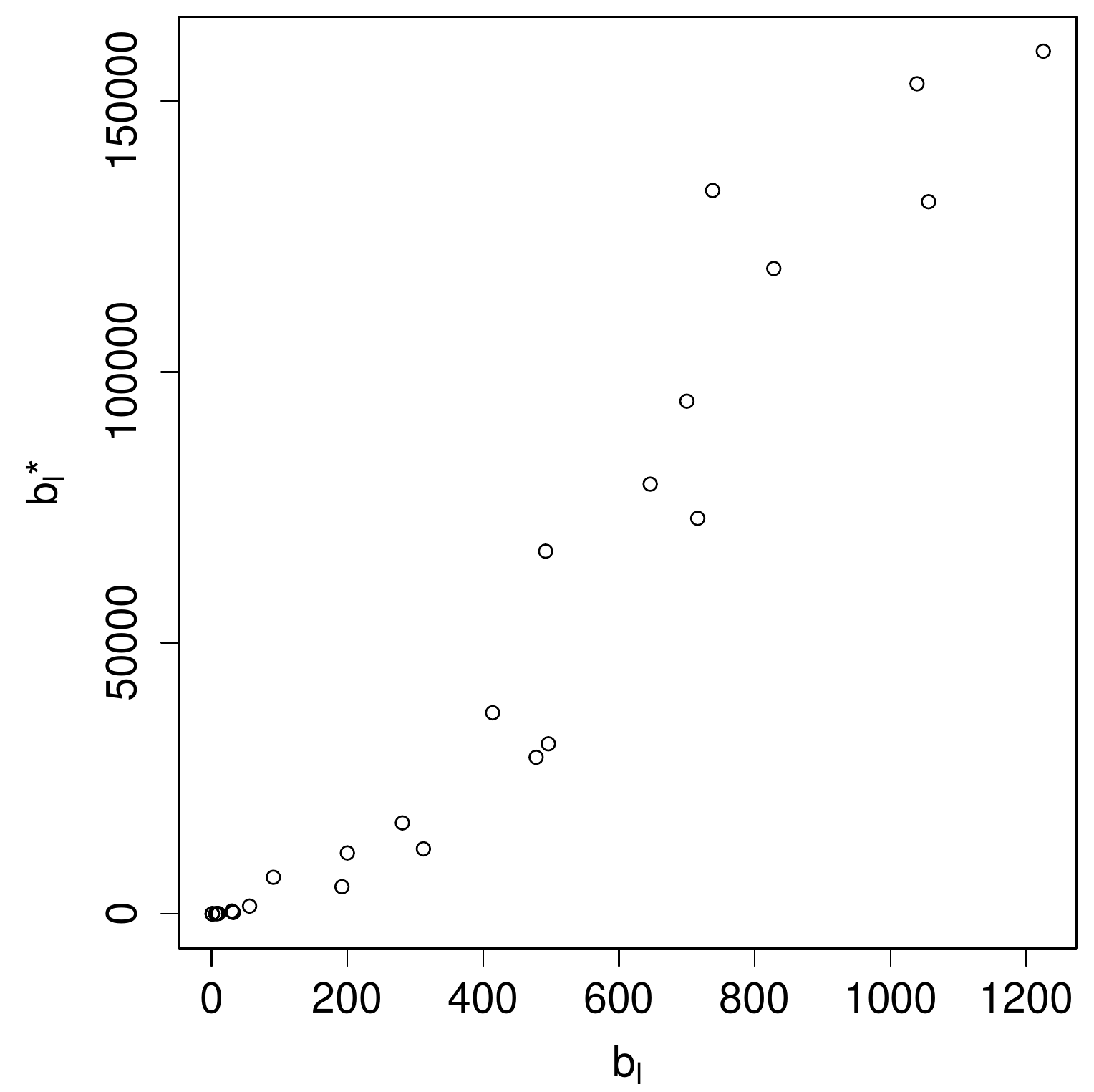}\hspace{0.5cm}\includegraphics[width=5.5cm]{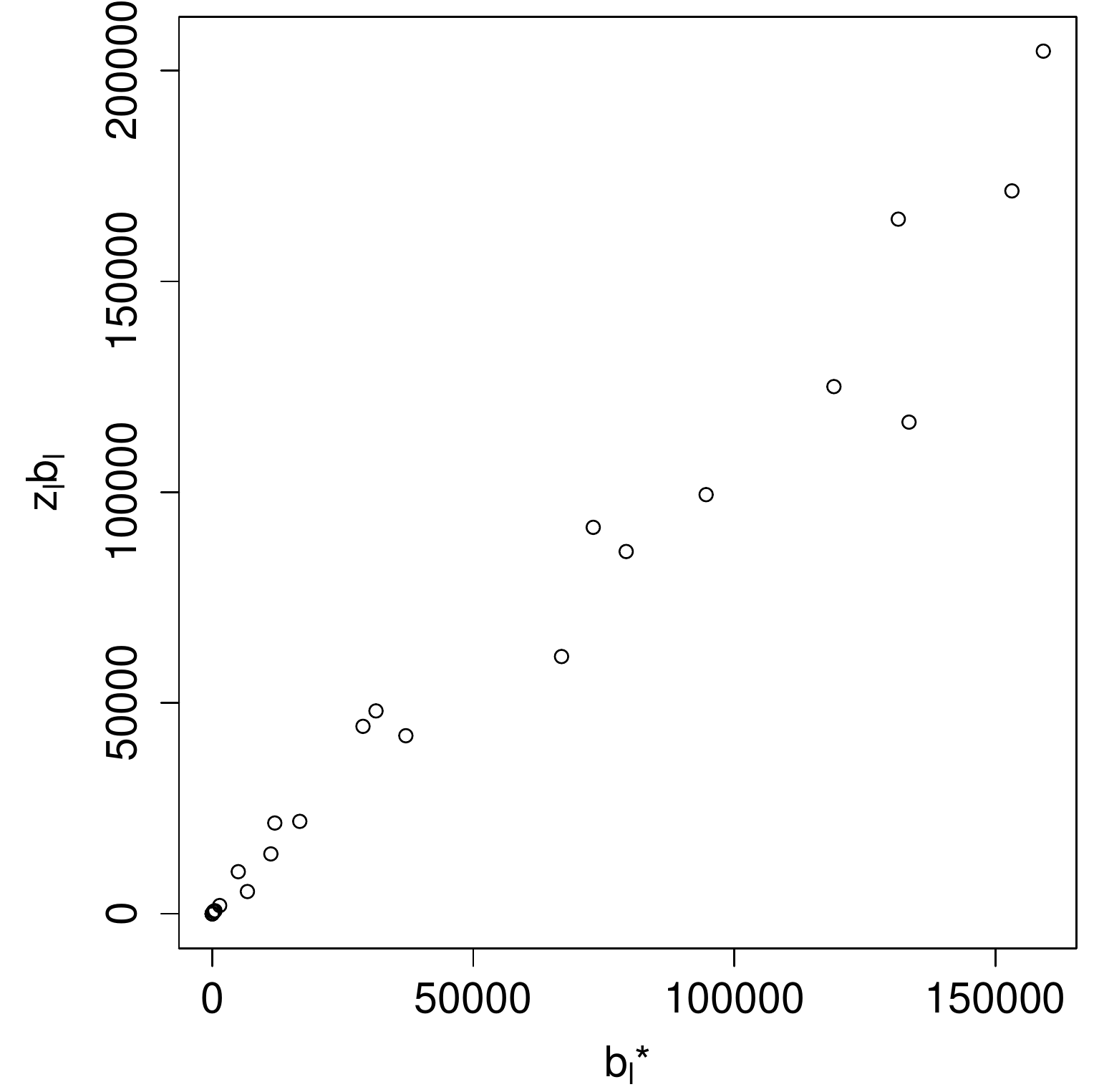}\\
  \caption{\label{complex} Evolution of $(b_l,b_l^*)$-left- and $(b_l^*,z_l b_l)$-right-, $l=0,\ldots,29$, given the samples $\overline{z}_{30}$ and $z_{30}$, and by considering $s_{max}=3$ and binomial control distributions.}
 \end{center}
\end{figure}

Hence, for each iteration of both methods we can determine the
order of the computational complexity as
$s_{max}\sum_{l=0}^{n-1}b_l$ and $s_{max}\sum_{l=0}^{n-1}b_l^*$,
respectively. Now, for each $l=0,1,\ldots,n-1$, $b_l$ depends on
the values of $s_{max}$, $\phi_l(Z_l)$ and $Z_{l+1}$, and $b_l^*$
on  $s_{max}$, $Z_l$ and $Z_{l+1}${, but it is not possible to
obtain closed forms of them}. To obtain an upper bound of $b_l$
one can obtain the dimension of the biggest { transition} tree. In
the case of binomial control functions, this tree can be generated
by considering $\phi_l(Z_l)=Z_l$ (the maximum number of
progenitors). Clearly, the dimension of this tree is greater than
or equal to that of the one that leads to $Z_{l+1}$. By an
empirical study (see supplementary material for details) we have
determined that $b_l=O(Z_l^{s_{max}-1})$. In a similar way, an
upper bound of $b_l^*$ is given by the dimension of the biggest
tree that can be generated by $Z_l$ individuals under the lack of
awareness of the exact number of progenitors $\phi_{l}(Z_l)$.
Again, we have determined empirically (see supplementary material
for details) that, in the case of binomial control functions,
$b_l^*=O(Z_l^{s_{max}})$. This fact allows us to compare, at least
roughly, the computational complexity of both methods, indicating
that for a generation of size $z$, one needs to generate trees of
dimension $z$ times { greater } when only the population size is
observed than when the number of progenitors, $\phi_l(z)$, is {
also available.  Figure \ref{complex}} reveals this fact in our
numerical example.
     This implies that the EM { procedure} requires much more time in each iteration when storing only $\mathcal{Z}_n$, compared to when storing $\overline{\mathcal{Z}}_n$. In particular, in our example, the same number of iterations of the procedure
   required a factor of 128 less time under the sample with
   observation of offspring and progenitors than under the sample
   based only on generation sizes.  Also, the
   second of these two procedures needed more iterations to
   reach convergence.  Hence, as was to be expected due to the
   relative loss of information, the second method is far
   more costly computationally than the first (by a
   factor of roughly 170 for a precision of
   $10^{-6}${, in terms of computational time}).  Moreover, this second procedure involves
   post--processing which involves running it several times for different seeds, and evaluating the exact likelihood at the convergence points.
\begin{remark}
  The example simulations were performed by parallel
   computing using the {\bf R} statistical software
   environment (see \cite{R}). For the
   estimator density and the
   exact log-likelihood function calculations, we used
   the \emph{sm} and \emph{polynom} packages (see \cite{sm}
   and \cite{polynom}), respectively.
\end{remark}

\section{Concluding Remarks}\label{sec:conclusions}

We have studied the maximum likelihood estimation of the
main parameters of the CBP with random control function
considering a nonparametric framework for the offspring
distribution and a parametric scheme for the control
process.  First, assuming the entire family tree is
observable, we determined the MLEs of the parameters
associated with the offspring distribution and with the control
law, and established their consistency and limiting
normality.  These results generalized those that had been
obtained for the parameters associated with the offspring law
for CBPs with deterministic control function.  We also
provided new results on the estimation of the control and
migration parameters, with particular note made of their
asymptotic properties.

Since in practice it is difficult to observe the entire family
tree, we considered two more realistic situations, one assuming
that the only observable data are the total number of individuals
and progenitors in each generation,
{$\overline{\mathcal{Z}}_n$}, and the other that even
only the generation sizes are observable, $\mathcal{Z}_n$. In both
cases, we addressed the problem of obtaining the MLEs of the main
parameters of the model by an incomplete data estimation
procedure. To this end, we made use of the EM algorithm.  A
simulated example showed that this seems to work appropriately
based on the sample {$\overline{\mathcal{Z}}_n$}. Based
on the sample {${\mathcal{Z}}_n$}, we encountered the
problem that the algorithm may converge to local maxima or saddle
points. In such a case, we proposed running the algorithm with a
large number of different starting values, and choosing the ones
associated with the highest value of the log-likelihood function
(this function can be evaluated although it can not be maximized
by standard methods). The simulated example showed this
methodological strategy to also work adequately. The procedure
based on knowledge of the total numbers of individuals and
progenitors converges rapidly, providing adequate accuracy with
reasonably short computation times. Storing only
{${\mathcal{Z}}_n$} however, we found the EM algorithm
to require not only much more time for each iteration but also
more iterations to reach convergence (with the same precision).

In the simulated example, we also illustrated the
consistency of the estimates based on the three samples.
(The only case established theoretically in the paper was
that corresponding to observing the entire family tree.)  We
then used a bootstrapping approach to get approximations to
the sampling distributions of the estimators obtained by the
EM algorithm, finding that the more information that the
samples contained, the smaller was the variability of the
estimator.

\section*{Acknowledgements}
{ The authors thank the referee for her/his careful
reading of our paper and for her/his constructive comments which
have improved its presentation. Also, the} authors would like to
thank Horacio Gonz\'alez-Velasco and Carlos Garc\'ia-Orellana for
providing them with computational support. Research supported by
the Ministerio de Econom\'ia y Competitividad and the FEDER
through the Plan Nacional de Investigaci\'on Cient\'ifica,
Desarrollo e Innovaci\'on Tecnol\'olgica, grant MTM2012-31235.

\section*{Appendix A. Proof of Theorem \ref{thm:MLE-complete}}\label{appendixA}
It is immediate to verify that the likelihood function based on $\mathcal{Z}_n^*$ is:
\begin{align*}
   \mathcal{L}( p,\theta \  | Z_l (k) &= z_l(k), 0\leq l \leq n-1; k\geq 0)=\\
&=\theta^{\sum_{l=0}^{n-1}\phi_l^*}A_{\sum_{l=0}^{n-1}z_l}(\theta)^{-1}\prod_{l=0}^{n-1} \frac{\phi_l^*! \ a_{z_l}(\phi_l^*)}{\prod_{k=0}^\infty z_l(k)!}\prod_{k=0}^\infty p_k^{z_l(k)},
\end{align*}
where $\phi_l^*=\sum_{k=0}^\infty z_l(k)$. Consequently, the log-likelihood function based on $\mathcal{Z}_n^*$ is:
\begin{equation}\label{eq:log-lik-tree}
\ell(p,\theta \ | \mathcal{Z}_n^*) =  f(p)+ g(\theta) + K ,
\end{equation}
with $f(p)=\sum_{l=0}^{n-1}\sum_{k=0}^\infty Z_l(k) \log
p_k$, $g(\theta)= \Delta_{n-1}\log \theta - \log
(A_{Y_{n-1}}(\theta))$ and $K$ some positive random variable
whose expression does not depend on $p$ or $\theta$.

From \eqref{eq:log-lik-tree}, one has to maximize $f(p)+g(\theta)$ subject to the constraints $\sum_{k=0}^\infty p_k=1$, $p_k\geq 0$, $k\geq 0$. Using the non-negativity of the Kullback-Leibler divergence, it can be verified that the value of $p$ which maximizes the function $f(p)\Delta_{n-1}$, and hence, $f(p)$, is

$$\widehat{p}_{k,n}= \frac{\sum_{l=0}^{n-1} Z_l(k)}{\sum_{l=0}^{n-1} \sum_{k=0}^\infty Z_l(k)}=\frac{Y_{n-1}(k)}{\Delta_{n-1}},\qquad k\geq 0.$$

Moreover, it can easily be shown that
$$\widehat{\theta}_n=\mu^{-1}\left(\frac{\Delta_{n-1}}{Y_{n-1}}\right)$$
is an extremum of the function $g(\theta)$. Taking into account that
$$\frac{d^2 g(\theta)}{d\theta^2}\Big|_{\theta=\widehat{\theta}_n}=-\frac{\Delta_{n-1}-\varepsilon(Y_{n-1},\theta)+\sigma^2(Y_{n-1},\theta)}{\theta^2}\Big|_{\theta=\widehat{\theta}_n}<0,$$
one has that $\widehat{\theta}_n$ is the maximum of $g(\theta)$ and then $(\widehat{p}_{n},\widehat{\theta}_n)$ maximizes $f(p)+g(\theta)$.


\section*{Appendix B. Proof of Theorem \ref{thm:asymp-prop-estim} }\label{appendixB}
\emph{(i)} We shall prove that $\widehat{p}_k$ is strongly
consistent for $p_k$, making use of a strong law of large
numbers for martingales. We shall fix $k\geq 0$, and prove that, as
$n\to\infty$,
\begin{equation}\label{eq:consistence-pk}
\widehat{p}_k=\frac{1}{\sum_{j=0}^{n-1}\phi_j(Z_j)}\sum_{i=1}^n \sum_{j=1}^{\phi_{i-1}(Z_{i-1})}I_{\{X_{i-1j}=k\}}\rightarrow p_k\quad a.s.\text{ on }\{Z_n\to\infty\}.
\end{equation}

For simplicity, we will consider $P[Z_n\to\infty]=1$. For each $i=1,2,\ldots$, let
\begin{eqnarray*}
V_i(k)&=&\sum_{j=1}^{\phi_{i-1}(Z_{i-1})} (I_{\{X_{i-1j}=k\}}-p_k),\\
\mathcal{H}_i&=&\sigma(X_{l-1j},\phi_{l-1}(k):1\leq l\leq i,j\geq 1,k\geq 0).
\end{eqnarray*}

It is verified that $\{V_i(k),\mathcal{H}_i\}_{i\geq 0}$ is a martingale difference. In these terms, $\widehat{p}_k-p_k=\Delta_{n-1}^{-1}\sum_{i=1}^n V_i(k)$.

For each $n\geq 0$, let $U_n=Y_{n-1}$. Then, taking into account Proposition \ref{prop:asymp-behaviour}(iv), to prove \eqref{eq:consistence-pk} one only needs to obtain that, as $n\to\infty$,
\begin{equation}\label{eq:lim-suma-vi}
U_n^{-1}\sum_{i=1}^n V_i(k)\rightarrow 0\quad a.s.
\end{equation}

Since $U_n\rightarrow\infty$, to prove \eqref{eq:lim-suma-vi}, using Theorem 2.18 in \cite{H-Hl}, it is enough to verify that $\sum_{i=1}^\infty U_i^{-2} E[|V_i(k)|^2| \mathcal{H}_{i-1}]<\infty\ a.s$. Now, let $M=\sup_{n\geq 0}\Delta_n Y_n^{-1}<\infty\ a.s.$ and $N=\sup_{n\geq 0}\varepsilon(Z_n)\phi_n(Z_n)^{-1}<\infty\ a.s.$ (guaranteed by Proposition \ref{prop:asymp-behaviour}(iv) and (vi), respectively). Then, one has

\begin{eqnarray*}
  \sum_{i=1}^\infty U_i^{-2} E[|V_i(k)|^2| \mathcal{H}_{i-1}] &=& \sum_{i=1}^\infty \frac{E\left[Var\left[\sum_{j=1}^{\phi_{i-1}(Z_{i-1})}I_{\{X_{i-1j}=k\}}\Big|\phi_{i-1}(Z_{i-1})\right]\right]}{Y_{i-1}^2} \\
   &=& \sum_{i=1}^\infty \frac{\varepsilon(Z_{i-1})p_k(1-p_k)}{Y_{i-1}^2} \\
   &=&  p_k(1-p_k)\sum_{i=1}^\infty \frac{\varepsilon(Z_{i-1})}{\phi_{i-1}(Z_{i-1})}\cdot\frac{\phi_{i-1}(Z_{i-1})}{\Delta_{i-1}^2}\cdot \left(\frac{\Delta_{i-1}}{Y_{i-1}}\right)^2\\
   &\leq&  p_k(1-p_k)N M^2 \sum_{i=1}^\infty \frac{1}{\phi_{i-1}(Z_{i-1})}<\infty\quad a.s.,
\end{eqnarray*}
where the last inequality is true due to $\phi_{i-1}(Z_{i-1})\leq\Delta_{i-1}$, $i\geq 1$, and Proposition \ref{prop:asymp-behaviour}(ii).

\vspace*{0.5cm}
The strong consistency of $\widehat{m}_n$ is a direct consequence of Proposition \ref{prop:asymp-behaviour}(ii)-(v).

\vspace*{0.2cm}
Taking into account that both $\widehat{m}_n$ and $\widehat{p}_k$ are strongly consistent for $m$ and $p_k$, respectively, on $\{Z_n\to\infty\}$, it is deduced that $\widehat{\sigma}_n^2$ is strongly consistent for $\sigma^2$.

\vspace*{0.2cm}
\emph{(ii)} The key to proving (ii)~(a) and (b) is to rewrite
$$  \widehat{p}_k - p_k   \stackrel{d}{=} \frac{1}{\Delta_{n-1}}\sum_{l=1}^{\Delta_{n-1}} \left(I_{\{X_{l}=k\}} - p_k\right),\quad
    \widehat{m}_n - m \stackrel{d}{=} \frac{1}{\Delta_{n-1}}\sum_{l=1}^{\Delta_{n-1}}(X_{l} - m),$$
    with $\stackrel{d}{=}$, as one recalls, denoting equal
    in distribution, and $\{X_l\}_{l\geq 1}$ being a
    sequence of i.i.d. random variables with common
    distribution being the offspring distribution. The results are
    derived by applying a central limit theorem for random
    sums as was done, \emph{mutatis mutandis}, in the
    proofs of Theorems 3.2 and 4.2 in \cite{art-2004b} for
    CBPs with deterministic control function.

Finally, to prove (ii)~(c), we  adapt the proof established in Theorem 3.1 in \cite{art-2005b} for CBPs with deterministic control function. We here provide just a brief scheme. The result is firstly proved for $\sum_{k=0}^\infty (k-m)^2\widehat{p}_k$, i.e., when $m$ is considered known. Due to the fact that one can write
\begin{eqnarray*}
\sum_{k=0}^\infty (k-m)^2 \widehat{p}_k
 &\stackrel{d}{=}& \frac{1}{\Delta_{n-1}}\sum_{l=1}^{\Delta_{n-1}} (X_{l}-m)^2,
\end{eqnarray*}
the result holds by using the central limit theorem cited
above, following similar steps to those in the proof of
Theorem 3.1 in \cite{art-2005b}. Now, notice that $
\sum_{k=0}^\infty (k-m)^2 \widehat{p}_k -
\widehat{\sigma}_n^2= (\widehat{m}_n - m)^2$, so that, by
considering (ii) (b),
$\sigma^2\Delta_{n-1}^{-1/2}\xrightarrow{P'} 0$, and
Slutsky's Theorem, one has $$\left(\sum_{k=0}^\infty (k-m)^2
  \widehat{p}_k -
  \widehat{\sigma}_n^2\right)\Delta_{n-1}^{1/2}
\xrightarrow{P'} 0, \mbox{ as } n\to\infty.$$ Hence, together with the fact that the result holds when $m$ is known and
Slutsky's Theorem, (ii)~(c) follows.

\section*{Appendix C. Proof of Theorem \ref{thm:asymp-prop-estim-theta-mu} }\label{appendixC}
\emph{(i)} This is immediate from Proposition \ref{prop:asymp-behaviour}(iv).

\emph{(ii)~(a)} For simplicity, we shall suppose $P[Z_n\to\infty]=1$. Let $D_i=\phi_{i-1}(Z_{i-1})-\mu(\theta)Z_{i-1}$ and
$\mathcal{F}_i=\sigma(Z_0,\ldots,Z_i,\phi_0(Z_0),\ldots,\phi_{i-1}(Z_{i-1}))$, $i=1,\ldots,n$, $n=1,2,\ldots$ We have
\begin{eqnarray}\label{eq:asympt-distrib-mu}
 Y_{n-1}^{1/2}\left(\widehat{\mu}-\mu(\theta)\right) &=&  Y_{n-1}^{-1/2} \Bigg[\sum_{i=1}^n \left((Z_{i-1}+1)^{1/2}-(\tau_m^{i-1}W)^{1/2}\right)\frac{D_i}{(Z_{i-1}+1)^{1/2}}\nonumber \\
&\phantom{=}&+W^{1/2}\sum_{i=1}^n \tau_m^{(i-1)/2} \frac{D_i}{(Z_{i-1}+1)^{1/2}}\Bigg]\nonumber
\end{eqnarray}
with $W$ being the limit variable introduced in \eqref{cond:prob-explosion}(c). Taking into account $\tau_m^{-n}Y_{n-1}\to (\tau_m -1)^{-1} W\ a.s.$ as $n\to\infty$, it follows that it is enough to prove
\begin{equation}\label{eq:t1-asympt-distrib-mu}
  \displaystyle(I)= \sum_{i=1}^n \left((Z_{i-1}+1)^{1/2}-(\tau_m^{i-1}W)^{1/2}\right)\frac{D_i}{(Z_{i-1}+1)^{1/2}}=o_P (\tau_m^{n/2})
\end{equation}
and
\begin{equation}\label{eq:t2-asympt-distrib-mu}
  \displaystyle(II)=(\tau_m -1)^{1/2}\sum_{i=1}^n \tau_m^{-(n-i+1)/2} \frac{D_i}{(Z_{i-1}+1)^{1/2}}\xrightarrow{d} N(0,\theta\mu'(\theta)),
\end{equation}
as $n\to\infty$, with $o_P(\cdot)$ denoting the stochastic order analogue of $o(\cdot)$ (i.e., write $X_n=o_P(Y_n)$ to mean $P(|X_n|\geq \epsilon |Y_n|)\to 0$, as $n\to\infty$, for each $\epsilon>0$). The proof follows similar steps to those given in \cite{Sriram}, Theorem 2.  For each $n\geq 0$, let
  $$A_n = \sum_{i=1}^n \tau_m^{(i-1)/2} \left(\left(\frac{Z_{i-1}+1}{\tau_m^{i-1}}\right)^{1/2}-W^{1/2}\right)^2\quad\mbox{and}\quad  B_n = \sum_{i=1}^n \tau_m^{(i-1)/2}\frac{D_i^2}{Z_{i-1}+1}.$$
Then, applying the Cauchy-Schwarz inequality, $|(I)|\leq A_n^{1/2} B_n^{1/2}$. By \eqref{cond:prob-explosion}(c), one obtains $(\tau_m^{-(i-1)}(Z_{i-1}+1))^{1/2}-W^{1/2}\to 0\ a.s.$, and consequently, using the Stolz-C\`{e}saro Lemma, $A_n=o_P\left(\sum_{i=1}^n \tau_m^{(i-1)/2}\right)=o_P\left(\tau_m^{n/2}\right)$. Now, using
\begin{equation}\label{eq:var-D-i}
    E[D_i^2|\mathcal{F}_{i-1}]=\theta \mu'(\theta) Z_{i-1},\quad i\geq 1,
\end{equation}
one has that $E[B_n]=O\left(\sum_{i=1}^n \tau_m^{(i-1)/2}\right)=O\left(\tau_m^{n/2}\right), \text{ as } n\to\infty,$
which implies that $|B_n|=O_P\left(\sum_{i=1}^n \tau_m^{(i-1)/2}\right)=O_P\left(\tau_m^{n/2}\right)$ as $n\to\infty$, with $O_P(\cdot)$ denoting the stochastic order analogue of $O(\cdot)$ (i.e., write $X_n=O_P(Y_n)$ to mean: for each $\epsilon>0$ there exists a real number $M$ such that $P(|X_n|\geq M |Y_n|)<\epsilon$ if $n$ is large enough). Hence \eqref{eq:t1-asympt-distrib-mu} follows.

To establish \eqref{eq:t2-asympt-distrib-mu}, let $\gamma_{nj}=D_{n-j+1}(Z_{n-j}+1)^{-1/2}$, $j=1,\ldots,n$. Then
\begin{align}\label{eq:t3-asympt-distrib-mu}
   (\tau_m -1)^{1/2}&\sum_{i=1}^n \tau_m^{-(n-i+1)/2} \frac{D_i}{(Z_{i-1}+1)^{1/2}} =  (\tau_m -1)^{1/2}\sum_{j=1}^n \tau_m^{-j/2} \frac{D_{n-j+1}}{(Z_{n-j}+1)^{1/2}}\nonumber \\
   &= (\tau_m -1)^{1/2}\left(\sum_{j=1}^J \tau_m^{-j/2} \gamma_{nj} + \sum_{j=J+1}^n \tau_m^{-j/2} \gamma_{nj}\right)\nonumber \\
   &= U_{Jn} + (\tau_m -1)^{1/2}\sum_{j=J+1}^n \tau_m^{-j/2} \gamma_{nj}=U_{nn},
\end{align}
with $U_{Jn}=(\tau_m -1)^{1/2}\sum_{j=1}^J \tau_m^{-j/2} \gamma_{nj}$, $J=1,\ldots,n$. For $J\geq 1$ and given $(t_1,\ldots,t_J)\in\R^J$, it can be proved, using analogous arguments to those given in the proof of Theorem 1 in \cite{Heyde-Brown-71}, jointly with the condition assumed  in \emph{(ii)}, that
$$E\left[e^{i \sum_{j=1}^J t_j \tau_m^{-j/2}\gamma_{nj}}\right]\rightarrow e^{-\frac{1}{2}\theta\mu'(\theta)\sum_{j=1}^J t_j^2 \tau_m^{-j}}, \text{ as }n\to\infty.$$
Consequently, for each $J=1,\ldots,n$, the vector $(\tau_m^{-1/2}\gamma_{n1},\ldots,\tau_m^{-J/2}\gamma_{nJ})$ is asymptotically multivariate normal as $n\to\infty$, and therefore $U_{Jn}\xrightarrow{d} U_J$, with $U_J$ following a $N(0,\theta\mu'(\theta)(\tau_m-1)\sum_{j=1}^J \tau_m^{-j})$. Now, from Chebyschev's inequality, \eqref{eq:var-D-i}, and \eqref{eq:t3-asympt-distrib-mu}, for every $n\geq 0$ and $\epsilon>0$, one has $P\left[|U_{Jn}-U_{nn}|>\epsilon\right] \leq \epsilon^{-2}(\tau_m-1)\theta\mu'(\theta)\sum_{j=J+1}^\infty \tau_m^{-j}$. In consequence, for some constant $k_0$,
$$\limsup_{n\to\infty} P\left[|U_{Jn}-U_{nn}|>\epsilon\right]\leq k_0 \sum_{j=J+1}^\infty \tau_m^{-j}\to 0, \text{ as } J\to\infty.$$

Therefore, from the fact that $U_J\xrightarrow{d} N(0,\theta\mu'(\theta))$ as $J\to\infty$, and Theorem 25.5 in \cite{Bi-79}, one obtains
$$U_{nn}\xrightarrow{d} N(0,\theta\mu'(\theta)),$$
as $n\to\infty$, and hence \eqref{eq:t2-asympt-distrib-mu} is proved.

\emph{(ii)~(b)} This is proved with identical arguments to those of \emph{(ii)~(a)}, setting in this case
$D_i=Z_i-\tau_m Z_{i-1}$ and $\mathcal{F}_i=\sigma(Z_0,\ldots,Z_i)$, $i=1,\ldots,n$, $n=1,2,\ldots$. Consequently,
$E[D_i^2|\mathcal{F}_{i-1}]=(\sigma^2\mu(\theta)+m^2\theta \mu '(\theta))Z_{i-1}$, $i\geq 1$, and now it is verified that
$$E\left[e^{i \sum_{j=1}^J t_j \tau_m^{-j/2}\gamma_{nj}}\right]\rightarrow e^{-\frac{1}{2}(\sigma^2\mu(\theta)+m^2\theta \mu '(\theta))\sum_{j=1}^J t_j^2 \tau_m^{-j}}, \text{ as }n\to\infty.$$

\section*{Supplementary material}
\section*{Simulated data}
We consider a CBP whose offspring distribution is given by
$p_0=0.1084$, $p_1=0.2709$, $p_2=0.3386$ and $p_3=0.2822$, and the
control variables $\phi_n(k)$ follow binomial distributions with
parameters $k$ and $q=0.6$. Thus, the offspring mean and variance
are $m=1.7946$ and $\sigma^2=0.9443$, respectively; $\theta=1.5$,
$\mu(\theta)=0.6$ and the mean growth rate is $\tau_m=1.0767$. We
simulated the first 30 generations of this process starting
with $Z_0=1$ individual. We denote by $z^*_{30}$,
$\overline{z}_{30}$ and $z_{30}$, the samples based on the entire
family tree, on the the total number of individuals and
progenitors in each generation, and on the generation sizes only,
respectively. The  data obtained are the following:

\begin{center}
\begin{longtable}{| c | c | c | c | c | c | c |}
  \hline
  $n$ & $Z_n$ & $\phi_n(Z_n)$ &$Z_n(0)$ & $Z_n(1)$ & $Z_n(2)$ & $Z_n(3)$  \\
  \hline
  0 & 1 & 1 & 0 & 1 & 0 & 0\\
1 & 1 & 1 & 0 & 1 & 0 & 0\\
2 & 1 & 1 & 0 & 0 & 1 & 0\\
3 & 2 & 1 & 0 & 0 & 1 & 0\\
4 & 2 & 2 & 0 & 0 & 1 & 1\\
5 & 5 & 5 & 0 & 2 & 2 & 1\\
6 & 9 & 6 & 1 & 2 & 2 & 1\\
7 & 9 & 7 & 2 & 2 & 1 & 2\\
8 & 10 & 8 & 0 & 3 & 1 & 4\\
9 & 17 & 14 & 1 & 8 & 3 & 2\\
10 & 20 & 14 & 0 & 8 & 2 & 4\\
11 & 24 & 17 & 2 & 2 & 6 & 7\\
12 & 35 & 25 & 1 & 6 & 8 & 10\\
13 & 52 & 39 & 4 & 11 & 14 & 10\\
14 & 69 & 48 & 10 & 15 & 13 & 10\\
15 & 71 & 38 & 8 & 9 & 14 & 7\\
16 & 58 & 36 & 1 & 5 & 17 & 13\\
17 & 78 & 51 & 5 & 13 & 15 & 18\\
18 & 97 & 61 & 7 & 28 & 13 & 13\\
19 & 93 & 61 & 5 & 22 & 22 & 12\\
20 & 102 & 64 & 6 & 15 & 20 & 23\\
21 & 124 & 72 & 7 & 14 & 25 & 26\\
22 & 142 & 76 & 5 & 24 & 32 & 15\\
23 & 133 & 73 & 9 & 21 & 22 & 21\\
24 & 128 & 81 & 9 & 16 & 33 & 23\\
25 & 151 & 86 & 8 & 21 & 34 & 23\\
26 & 158 & 83 & 7 & 19 & 34 & 23\\
27 & 156 & 94 & 11 & 26 & 32 & 25\\
28 & 165 & 94 & 11 & 29 & 24 & 30\\
29 & 167 & 107 & 10 & 27 & 37 & 33\\
30 & 200 & $\cdot$   &  $\cdot$ &   $\cdot$ & $\cdot$  & $\cdot$\\
\cline{2-2}
&$z_{30}$& &&&&\\\cline{2-3}
&\multicolumn{2}{c|}{$\overline{z}_{30}$}& &&&\\\cline{2-7}
\multicolumn{3}{|c|}{}&\multicolumn{4}{c|}{$z^*_{30}$} \\
\hline\caption{Simulated data}\label{sd}
\end{longtable}
\end{center}
\section*{Analysis of the robustness of the EM algorithm based on the sample $z_{30}$}

The EM algorithm based on the sample $z_{30}$ is observed not to
be at all robust to the choice of the initial values
$(p^{(0)},\theta^{(0)})$, with convergence to different estimates
that could be local maxima or saddle points of the log-likelihood
function. We detected this fact by starting the algorithm with 300
different initial values, choosing each $p^{(0)}$ randomly from a
Dirichlet distribution with all the parameters equal to one (that
is, uniformly from the unit simplex) and each $\theta^{(0)}$
through the equation $\theta^{(0)}=q^{(0)}(1-q^{(0)})^{-1}$, with
$q^{(0)}$ sampled from a uniform distribution on the open interval
$(0,1)$. To overcome this problem, we propose in the paper a methodological approach in
order to choose the best approximation to  the MLE based on
$\mathcal{Z}_{n}$ -called EM estimate-, that  we analyze in
detail below.

The log-likelihood function based on the sample $\mathcal{Z}_n=\{Z_0,\ldots,Z_n\}$, denoted by
 $\ell(p,\theta\mid \mathcal{Z}_n)$, is given, in the case of binomial control functions,  by
\begin{equation}\label{logexa}
    \ell(p,\theta\mid Z_l=z_l, \ l=0,\ldots, n)=\sum_{j=0}^{n-1}\log\left(\sum_{l=0}^{z_j}P_{z_{j+1}}^{*l}\binom{z_j}{l}
    \frac{\theta^l}{(1+\theta)^{z_j}}\right)
\end{equation}
with $P_{\cdot}^{*l}$ denoting the   $l$-fold convolution of the
offspring law $p$. We assume, for computational purposes, that the
support of this distribution is $\{0,\ldots,s_{max}\}$, with
$s_{max}$ denoting the maximum number of offspring per progenitor
(in our example $s_{max}=3$). Notice that $P_{z_{j+1}}^{*l}$ is
the coefficient of the monomial of degree $z_{j+1}$ of the
polynomial $(\sum_{k=0}^{s_{max}}p_ks^k)^l$. This fact allows us to
develop a computational program to evaluate the log--likelihood
function on each pair $(p,\theta)$.

\begin{figure}[h]
\begin{center}
\includegraphics*[width=5cm]{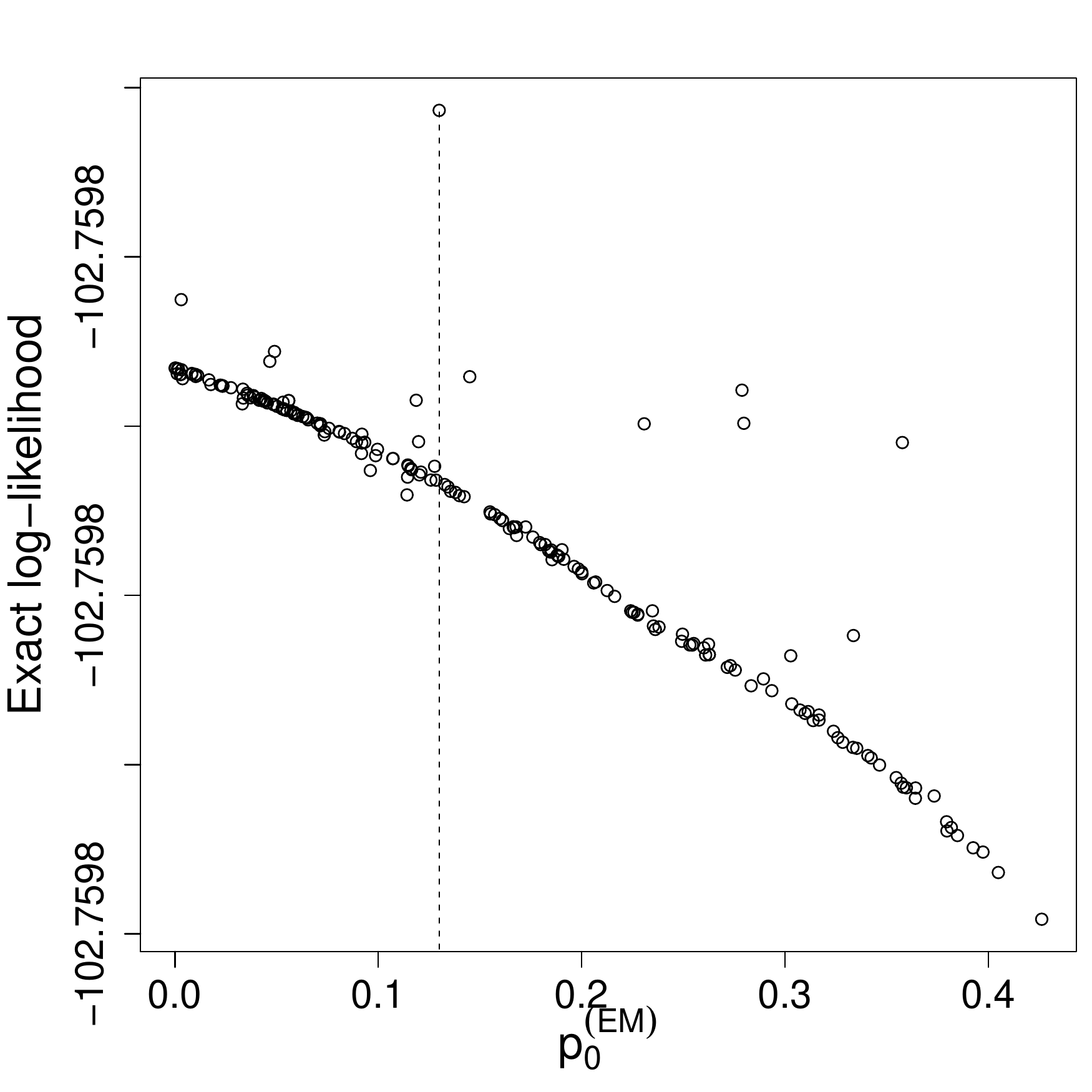}\includegraphics*[width=5cm]{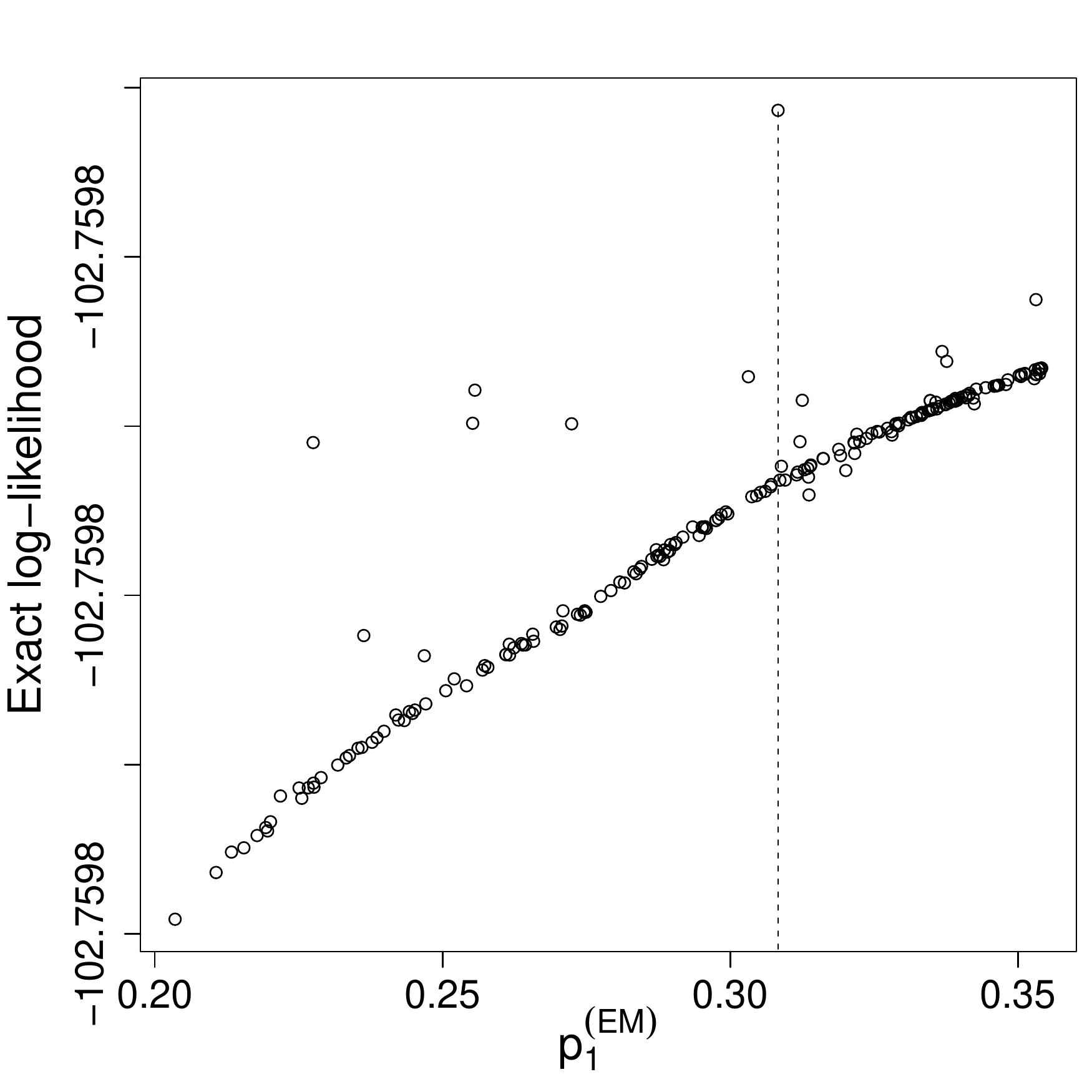}\includegraphics*[width=5cm]{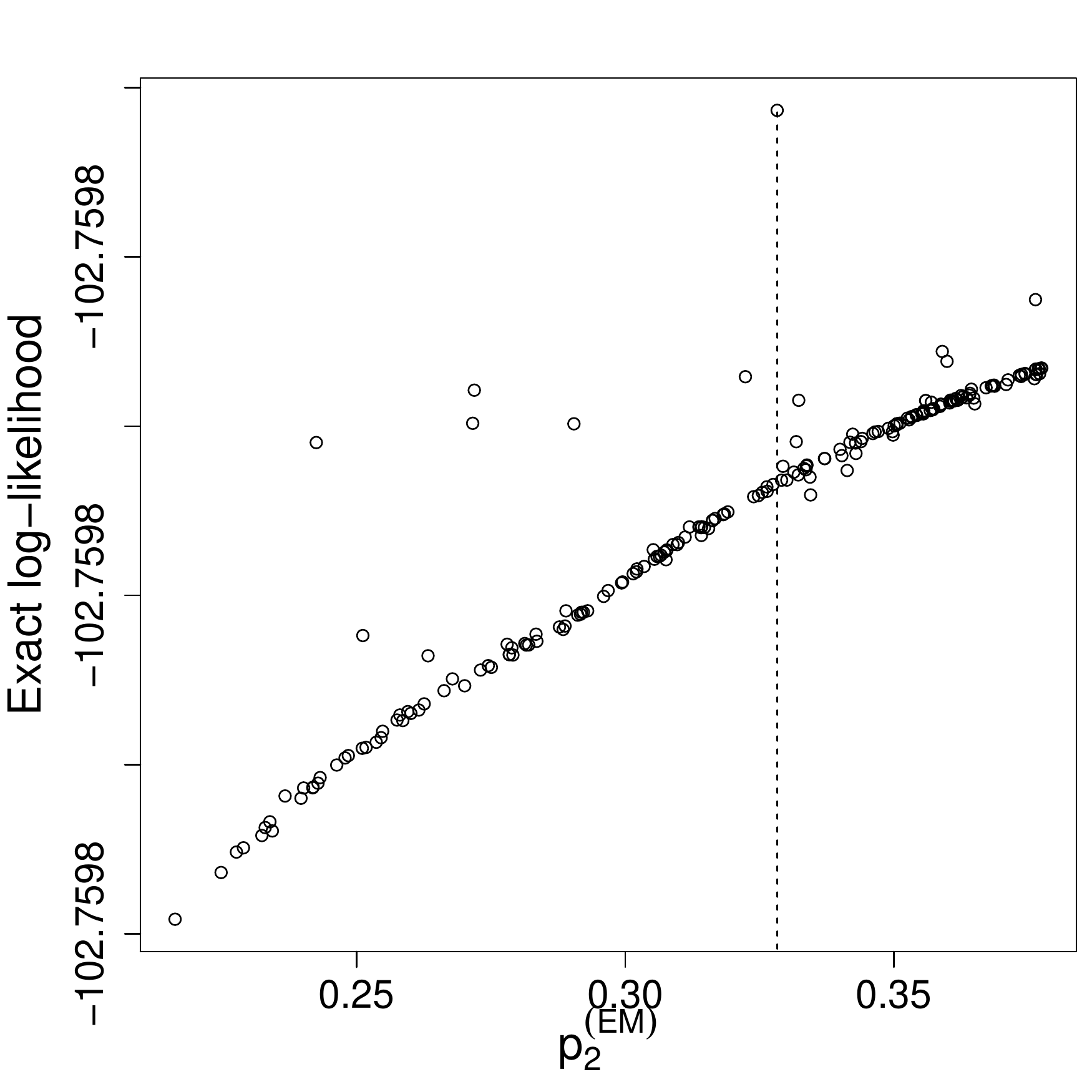}
\caption{Exact log-likelihood function versus the convergence
points of the EM algorithm for the parameters $p_0$, $p_1$ and
$p_2$, denoted by $p_0^{(EM)}$, $p_1^{(EM)}$ and
$p_2^{(EM)}$.\label{exlike} }
\end{center}
\end{figure}

\begin{figure}[h]
\begin{center}
\includegraphics*[width=5cm]{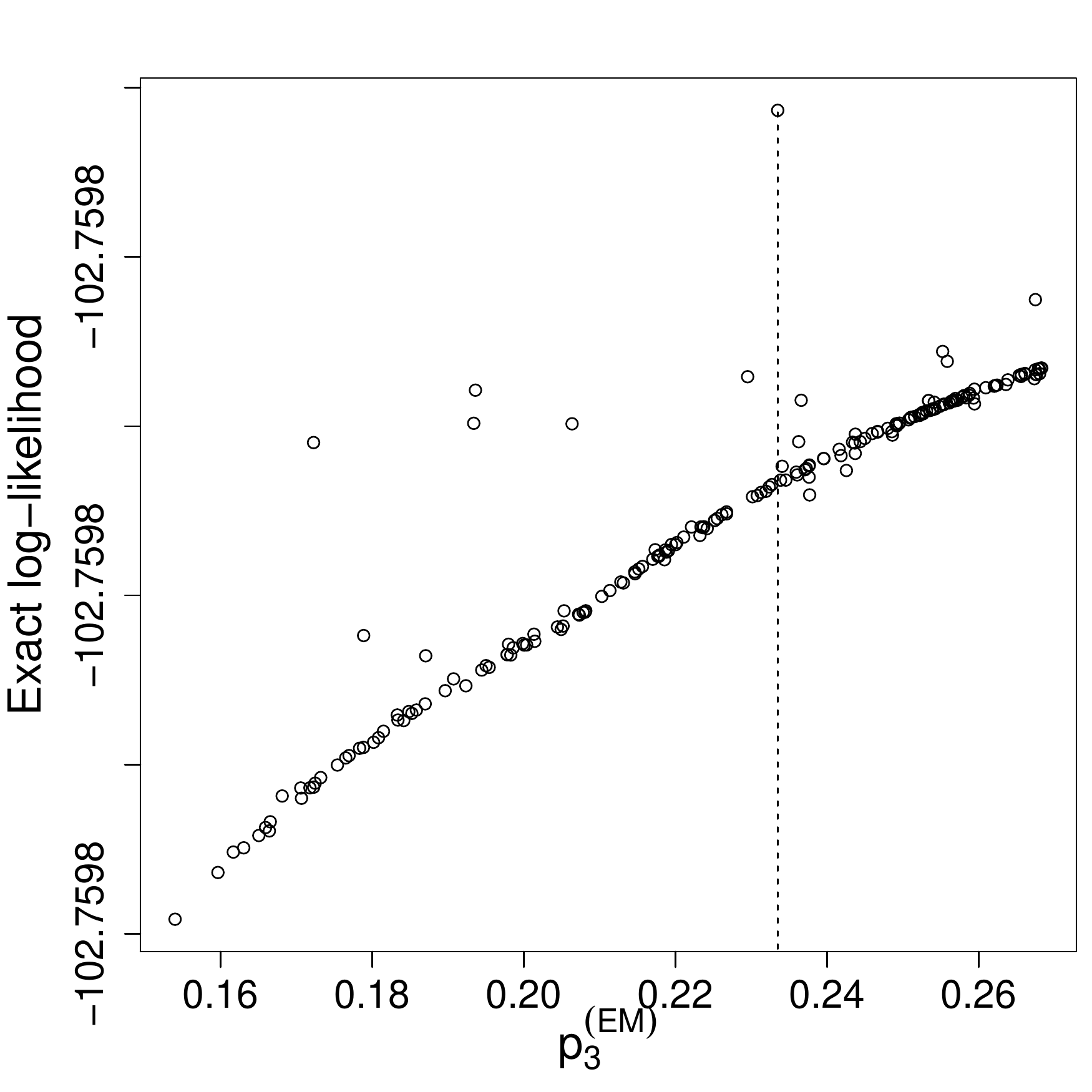}\hspace*{1cm}\includegraphics*[width=5cm]{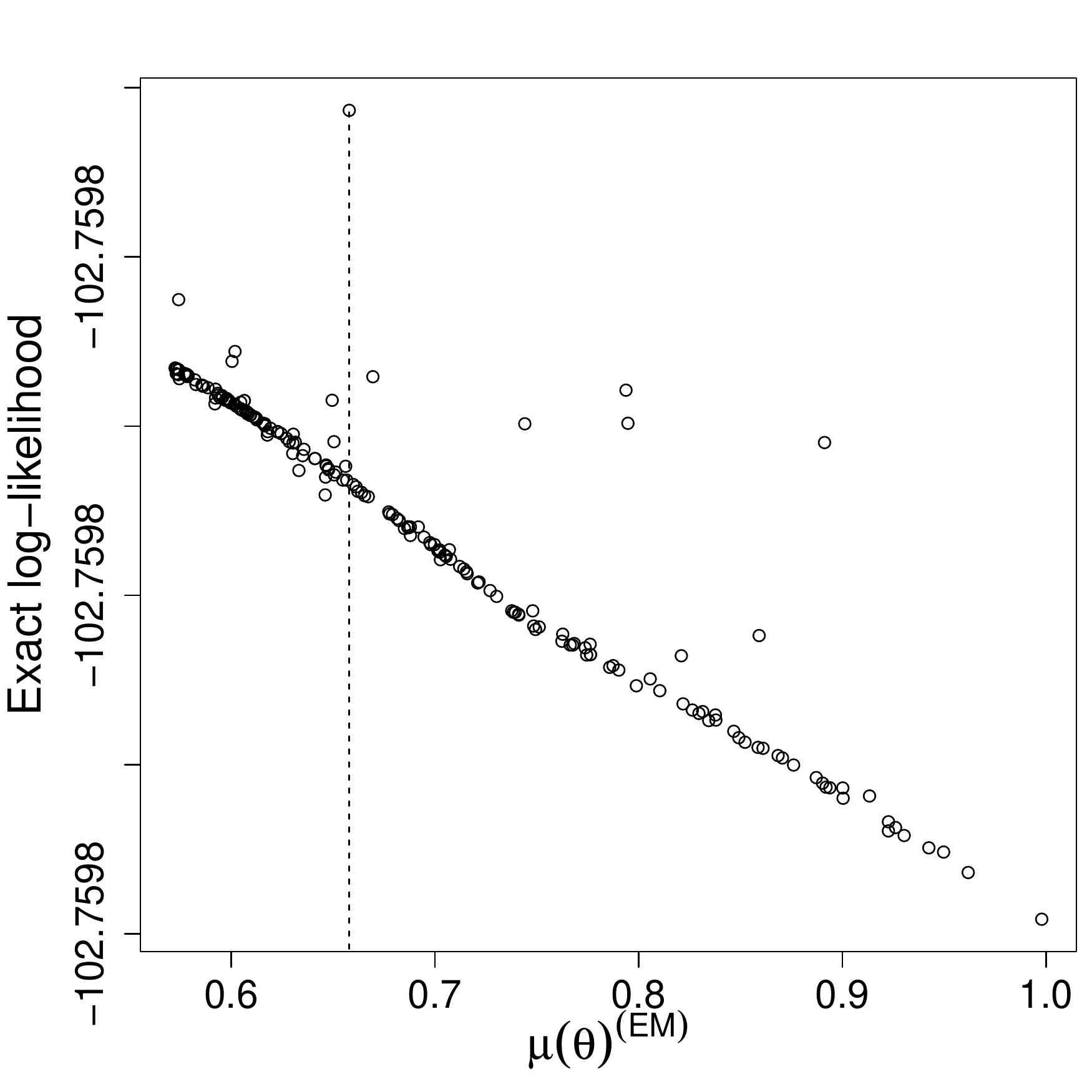}
\caption{Exact log-likelihood function versus the convergence
points of the EM algorithm for the parameters $p_3$ and
$\mu(\theta)$, denoted by $p_3^{(EM)}$ and
$\mu(\theta)^{(EM)}$.\label{exlike1} }
\end{center}
\end{figure}

Thus, the methodological approach we propose consists of taking
as EM estimates of the parameters those associated with the
greatest log-likelihood when is evaluated  at the convergence
points of the EM algorithm started with different randomly chosen
values of the parameters. Related to our example, in the following
figures we show the exact values of the  log-likelihood function
versus the convergence points of the EM algorithm started with
each one of the 300 different seeds, for the parameters $p_0$
(Figure \ref{exlike}, left), $p_1$ (Figure \ref{exlike}, center),
$p_2$ (Figure \ref{exlike}, right), $p_3$ (Figure \ref{exlike1},
left) and $\mu(\theta)$ (Figure \ref{exlike1}, right).

The values that maximizes the  log-likelihood function (\ref{logexa})
are given in Table \ref{t1}, and shown in Figures \ref{exlike} and
\ref{exlike1} with vertical dashed lines.
\begin{table}
  \centering
  \begin{tabular}{ccccc||cc||cc}
   \cline{2-9}
\multicolumn{1}{c}{}&\multicolumn{7}{c}{PARAMETERS}\\
\hline
       SAMPLE & $p_0$ & $p_1$ & $p_2$ & $p_3$ & $m$ & $\sigma^2$ & $\mu(\theta)$ & $\tau_m$ \\\hline
    $z_{30}$ & .1299 & .3083 & .3283 & .2335 & 1.6653 & .9496 & .6579 & 1.0957 \\\hline
TRUE VALUE     & .1084 & .2709 & .3386 & .2822 & 1.7946 & .9443 & .6000 & 1.0767 \\
      \end{tabular}
      \caption{Estimates of the parameters of interest based on the sample $z_{30}$.}\label{t1}
\label{t1}
\end{table}

\begin{figure}[h]
\begin{center}
\includegraphics*[width=6cm]{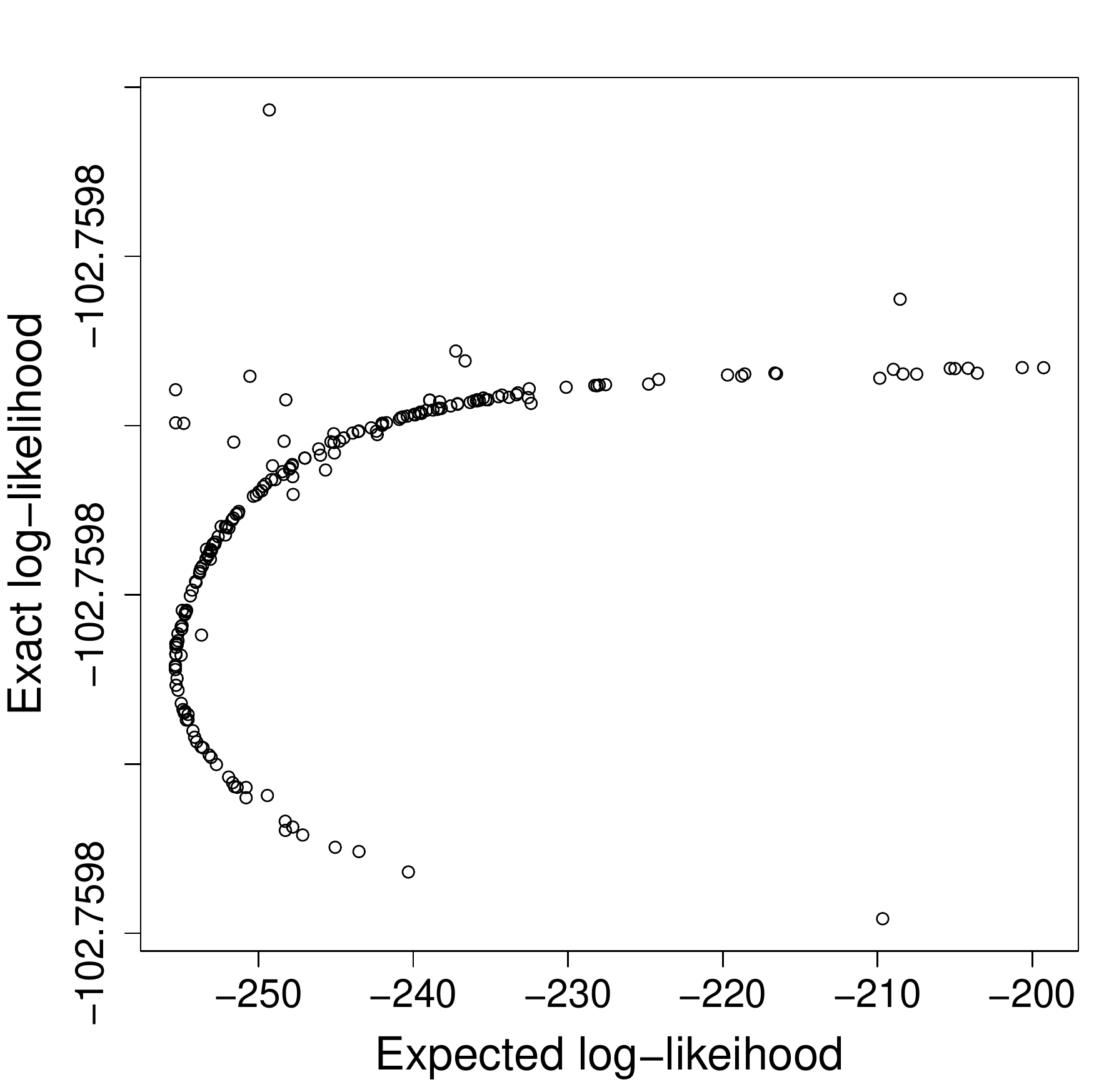}
\caption{Exact log-likelihood function versus   expected log-likelihood.\label{exexac} }
\end{center}
\end{figure}

Finally, it is also
worth mentioning that the expectation of the log-likelihood
$\ell(p,\theta \ | \mathcal{Z}_n^*,\ {\mathcal{Z}}_n)$ with
respect to the distribution
 $\mathcal{Z}_n^*|({\mathcal{Z}}_n,\{p^{(EM)},\theta^{(EM)}\})$, with $(p^{(EM)},$ $\theta^{(EM)})$ a convergence point of the EM algorithm --this kind of expected values can be calculated in each iteration of the algorithm (see Equation (7) in the paper)--,  can not be used to determine the maximum likelihood estimates, as an alternative to our proposal. This is due to the fact that it can happen that $$E_{\mathcal{Z}_n^*\mid (\mathcal{Z}_n,\{ p^{(EM)},\theta^{(EM)}\})}[\ell(p,\theta\mid \mathcal{Z}_n^*, \mathcal{Z}_n)]\geq E_{\mathcal{Z}_n^*\mid (\mathcal{Z}_n,\{ \tilde{p}^{(EM)},\tilde{\theta}^{(EM)}\})}[\ell(p,\theta\mid \mathcal{Z}_n^*, \mathcal{Z}_n)]$$ and $$\ell(p^{(EM)},\theta^{(EM)}\mid \mathcal{Z}_n)< \ell(\tilde{p}^{(EM)},\tilde{\theta}^{(EM)}\mid \mathcal{Z}_n),$$ being $(p^{(EM)},\ \theta^{(EM)})$ and  $(\tilde{p}^{(EM)},$ \ $\tilde{\theta}^{(EM)})$ two convergence points provided by the EM algorithm.
 Figure \ref{exexac} shows this fact. We plot on it $\ell(p^{(EM)},\theta^{(EM)} \ |
 z_{30})$ versus $E_{\mathcal{Z}_{30}^*\mid (z_{30},\{ p^{(EM)},\theta^{(EM)}\})}[\ell(p,\theta\mid \mathcal{Z}_{30}^*,
z_{30})]$, with $(p^{(EM)},\theta^{(EM)})$ the convergence points
of the EM algorithm started with the 300  randomly chosen seeds.

\section*{Computational complexity}
Focussing our interest in the case of binomial control functions, in order to determine  upper bounds of the values $b_l$ and $b_l^*$, let us introduce the following functions.  Let $b(z_l,\phi_l^*, z_{l+1},s_{max})$ the function that provides the number of possible vectors $(z_l(0),\ldots, z_l(s_{max}))$  such that $\sum_{k=0}^{s_{max}}z_l(k)=\phi_l^*$ and
$\sum_{k=0}^{s_{max}}kz_l(k)=z_{l+1}$, with  $z_l$, $\phi_l^*$ and $z_{l+1}$ whatever  possible values of the variables $Z_l$, $\phi_l(Z_l)$ and $Z_{l+1}$, respectively, and $s_{max}$  with the maximum
number of offspring per progenitor. Notice that given the sample $\overline{\mathcal{Z}}_n$, then $b_l=b(Z_l,\phi_l(Z_l), Z_{l+1},s_{max})$, $l=0,1,\ldots, n-1$. It is also remarkable that $b_l$ depends on $Z_l$ through $\phi_l(Z_l)$. Analogously, let  $b^*(z_l, z_{l+1},s_{max})$ the function that provides the number of possible vectors $(z_l(0),\ldots, z_l(s_{max}))$  such that $1\leq \sum_{k=0}^{s_{max}}z_l(k)\leq z_l$ and
$\sum_{k=0}^{s_{max}}kz_l(k)=z_{l+1}$ , with  $z_l$, and $z_{l+1}$ whatever possible values of the variables $Z_l$ and $Z_{l+1}$, respectively (we have not considered the case $\sum_{k=0}^{s_{max}}z_l(k)=0$ because it does not lead to any factible value in the case $z_{l+1}=0$ or to the null vector if $z_{l+1}=0$; in any case its contribution is irrelevant for our purpose). Notice that $b^*(z_l, z_{l+1},s_{max})=\sum_{1\leq \phi_l^*\leq z_l}b(z_l,\phi_l^*, z_{l+1},s_{max})$. Assuming the sample ${\mathcal{Z}}_n$, $b_l^*=b^*(Z_l, Z_{l+1},s_{max})$, $l=0,1,\ldots, n-1$.

Thus, to obtain  upper bounds of the functions $b$ and $b^*$ in terms of $z_l$ and $s_{max}$, we determined the functions
\begin{eqnarray}
    b_{max}&=&b_{max}(z_l,s_{max})=\max_{\tiny{\begin{array}{c}1\leq \phi_l^*\leq z_l\\ 0\leq z_{l+1}\leq s_{max}\cdot\phi_l^*\end{array}}}b(z_l,\phi_l^*,z_{l+1},s_{max})\nonumber\\&=&
    \max_{ 0\leq z_{l+1}\leq s_{max}\cdot z_l}b(z_l,z_l,z_{l+1},s_{max})\label{nn1}
\end{eqnarray}
and
\begin{equation}\label{nn2}
    b_{max}^*=b_{max}^*(z_l,s_{max})=
    \max_{ 0\leq z_{l+1}\leq s_{max} \cdot z_l}b^*(z_l,z_{l+1},s_{max}).
\end{equation}
To get these maximum values, we have considered three possible values of $s_{max}=3,4,5$, and for each one we have obtained the values of the function $b$ for $\phi_l^*$ going from 1 to 167 (this is the maximum value of $\phi_l(Z_l)$ in our simulated sample -see Table \ref{sd}) and for  $z_{l+1}$ going from 0 to $167\cdot s_{max}$. The values obtained have been stored in matrices of dimension $(167\cdot s_{max}+1)\times 167$. Each column corresponds to a possible value of $\phi_l^*$, going from 1 to 167, and each row to one of $z_{l+1}$, from 0 to $167\cdot s_{max}$. Notice that the non-null values for the column corresponding to $\phi_l^*$ are the $\phi_l^*\cdot s_{max}+1$ first elements. The matrices obtained are stored in the files \texttt{tree-max-3.cvs} (for $s_{max}=3$),  \texttt{tree-max-4.cvs} (for $s_{max}=4$) and \texttt{tree-max-5.cvs} (for $s_{max}=5$). From them, and taking into account (\ref{nn1}) and (\ref{nn2}), it is easy to obtain the values of  $b_{max}(z_l,s_{max})$ and $b_{max}^*(z_l,s_{max})$, which are given in Table \ref{t2}.  For each value of $s_{max}$, analysing the relationship between $z_l$ and $b_{max}(z_l,s_{max})$ and $z_l$ and $b_{max}^*(z_l,s_{max})$, using polynomial regression methods, it can be concluded that $b_{max}(z_l,s_{max})=O(z_l^{s_{max}-1})$ and  $b_{max}^*(z_l,s_{max})=O(z_l^{s_{max}})$. Consequently, we infer that  $b_l=O(Z_l^{s_{max}-1})$ and  $b_l^*=O(Z_l^{s_{max}})$.
\pagebreak
\begin{center}
{\tiny
\begin{longtable}{ r || r | r || r | r || r | r }
 & \multicolumn{2}{c||}{$s_{max}=3$}&
\multicolumn{2}{c||}{$s_{max}=4$}&
\multicolumn{2}{c}{$s_{max}=5$}\\
\hline
$z_l$&$b_{max}$&$b_{max}^*$&$b_{max}$&$b_{max}^*$&$b_{max}$&$b_{max}^*$\\
\hline
1&1&1&1&1&1&1\\
2&2&3&3&4&3&4\\
3&3&6&5&8&6&9\\
4&5&9&8&14&12&19\\
5&6&15&12&24&20&36\\
6&8&22&18&37&32&63\\
7&10&29&24&58&49&103\\
8&13&39&33&85&73&164\\
9&15&51&43&117&102&249\\
10&18&66&55&164&141&369\\
11&21&84&69&218&190&525\\
12&25&103&86&287&252&736\\
13&28&124&104&372&325&1006\\
14&32&150&126&473&414&1355\\
15&36&178&150&598&521&1790\\
16&41&213&177&736&649&2332\\
17&45&249&207&906&795&3000\\
18&50&286&241&1102&967&3809\\
19&55&331&277&1326&1165&4789\\
20&61&378&318&1585&1394&5953\\
21&66&433&362&1875&1651&7337\\
22&72&492&410&2210&1944&8965\\
23&78&552&462&2586&2275&10873\\
24&85&618&519&3002&2649&13091\\
25&91&691&579&3478&3061&15653\\
26&98&769&645&3997&3523&18603\\
27&105&856&715&4575&4035&21982\\
28&113&945&790&5217&4604&25833\\
29&120&1036&870&5923&5225&30213\\
30&128&1140&956&6706&5910&35153\\
31&136&1246&1046&7545&6660&40728\\
32&145&1366&1143&8475&7483&46986\\
33&153&1489&1245&9486&8372&54003\\
34&162&1614&1353&10585&9343&61824\\
35&171&1750&1467&11779&10395&70533\\
36&181&1893&1588&13062&11538&80195\\
37&190&2046&1714&14456&12764&90880\\
38&200&2209&1848&15956&14090&102681\\
39&210&2374&1988&17565&15516&115675\\
40&221&2545&2135&19309&17053&129965\\
41&231&2731&2289&21161&18691&145621\\
42&242&2921&2451&23146&20451&162758\\
43&253&3129&2619&25271&22330&181469\\
44&265&3340&2796&27544&24342&201853\\
45&276&3553&2980&29976&26476&224027\\
46&288&3784&3172&32537&28754&248116\\
47&300&4021&3372&35277&31174&274220\\
48&313&4274&3581&38181&33751&302490\\
49&325&4536&3797&41269&36471&333023\\
50&338&4801&4023&44542&39361&365983\\
51&351&5077&4257&47991&42416&401493\\
52&365&5368&4500&51647&45654&439716\\
53&378&5668&4752&55500&49060&480793\\
54&392&5987&5014&59569&52662&524907\\
55&406&6309&5284&63877&56455&572201\\
56&421&6634&5565&68388&60459&622833\\
57&435&6985&5855&73135&64656&676982\\
58&450&7339&6155&78124&69079&734837\\
59&465&7717&6465&83383&73720&796574\\
60&481&8102&6786&88913&78602&862397\\
61&496&8490&7116&94674&83705&932496\\
62&512&8896&7458&100732&89064&1007118\\
63&528&9316&7810&107061&94671&1086429\\
64&545&9751&8173&113710&100551&1170652\\
65&561&10204&8547&120663&106681&1260022\\
66&578&10661&8933&127909&113101&1354759\\
67&595&11126&9329&135487&119799&1455114\\
68&613&11617&9738&143377&126804&1561303\\
69&630&12113&10158&151630&134091&1673615\\
70&648&12640&10590&160256&141702&1792281\\
71&666&13171&11034&169209&149625&1917588\\
72&685&13706&11491&178531&157891&2049806\\
73&703&14267&11959&188219&166471&2189195\\
74&722&14839&12441&198341&175413&2336068\\
75&741&15434&12935&208875&184701&2490790\\
76&761&16045&13442&219772&194370&2653595\\
77&780&16660&13962&231113&204389&2824802\\
78&800&17290&14496&242854&214808&3004715\\
79&820&17945&15042&255093&225610&3193697\\
80&841&18611&15603&267785&236833&3392022\\
81&861&19307&16177&280916&248442&3600097\\
82&882&20008&16765&294534&260493&3818240\\
83&903&20713&17367&308599&272965&4046845\\
84&925&21454&17984&323218&285900&4286352\\
85&946&22202&18614&338373&299260&4537044\\
86&968&22982&19260&354009&313104&4799346\\
87&990&23774&19920&370176&327410&5073587\\
88&1013&24571&20595&386859&342223&5360240\\
89&1035&25391&21285&404172&357501&5659681\\
90&1058&26233&21991&422074&373309&5972353\\
91&1081&27094&22711&440500&389620&6298680\\
92&1105&27983&23448&459547&406484&6639100\\
93&1128&28877&24200&479155&423856&6994059\\
94&1152&29780&24968&499457&441804&7364046\\
95&1176&30722&25752&520404&460300&7749502\\
96&1201&31669&26553&541966&479397&8150987\\
97&1225&32659&27369&564199&499045&8569019\\
98&1250&33656&28203&587040&519319&9004036\\
99&1275&34658&29053&610656&540186&9456496\\
100&1301&35693&29920&635009&561704&9926962\\
101&1326&36746&30804&660026&583820&10415946\\
102&1352&37827&31706&685766&606612&10924007\\
103&1378&38932&32624&712211&630045&11451680\\
104&1405&40043&33561&739484&654181&11999519\\
105&1431&41171&34515&767575&678962&12568160\\
106&1458&42335&35487&796381&704473&13158290\\
107&1485&43511&36477&826014&730674&13770375\\
108&1513&44730&37486&856405&757632&14405000\\
109&1540&45955&38512&887694&785285&15062801\\
110&1568&47186&39558&919870&813722&15744388\\
111&1596&48461&40622&952866&842901&16450391\\
112&1625&49748&41705&986747&872893&17181481\\
113&1653&51074&42807&1021440&903631&17938314\\
114&1682&52419&43929&1057121&935211&18721555\\
115&1711&53770&45069&1093798&967585&19531923\\
116&1741&55148&46230&1131352&1000830&20370041\\
117&1770&56558&47410&1169850&1034875&21236691\\
118&1800&57989&48610&1209272&1069820&22132755\\
119&1830&59459&49830&1249738&1105615&23058898\\
120&1861&60936&51071&1291297&1142341&24015730\\
121&1891&62421&52331&1333791&1179921&25004058\\
122&1922&63959&53613&1377348&1218463&26024628\\
123&1953&65504&54915&1421889&1257911&27078204\\
124&1985&67100&56238&1467550&1298352&28165598\\
125&2016&68708&57582&1514386&1339705&29287641\\
126&2048&70323&58948&1562277&1382082&30445101\\
127&2080&71976&60334&1611296&1425424&31638919\\
128&2113&73655&61743&1661361&1469823&32870171\\
129&2145&75366&63173&1712643&1515192&34139421\\
130&2178&77111&64625&1765228&1561651&35447522\\
131&2211&78863&66099&1818932&1609135&36795378\\
132&2245&80633&67596&1873831&1657742&38183913\\
133&2278&82452&69114&1929902&1707380&39613984\\
134&2312&84282&70656&1987249&1758174&41086588\\
135&2346&86169&72220&2046012&1810056&42602574\\
136&2381&88064&73807&2105960&1863129&44162971\\
137&2415&89966&75417&2167236&1917295&45768701\\
138&2450&91919&77051&2229752&1972687&47420788\\
139&2485&93891&78707&2293624&2029230&49120185\\
140&2521&95907&80388&2359009&2087034&50868375\\
141&2556&97950&82092&2425713&2145996&52665971\\
142&2592&100001&83820&2493818&2206254&54513965\\
143&2628&102081&85572&2563232&2267730&56413420\\
144&2665&104204&87349&2634106&2330539&58365425\\
145&2701&106349&89149&2706640&2394571&60371037\\
146&2738&108546&90975&2780565&2459973&62431379\\
147&2775&110751&92825&2855965&2526660&64547499\\
148&2813&112964&94700&2932815&2594754&66720601\\
149&2850&115242&96600&3011200&2664140&68951927\\
150&2888&117530&98526&3091346&2734970&71242921\\
151&2926&119876&100476&3172971&2807155&73594370\\
152&2965&122241&102453&3256219&2880823&76007435\\
153&3003&124614&104455&3340992&2955852&78483347\\
154&3042&127029&106483&3427378&3032403&81023258\\
155&3081&129480&108537&3515649&3110380&83628483\\
156&3121&131965&110618&3605542&3189918&86300210\\
157&3160&134495&112724&3697138&3270889&89039747\\
158&3200&137034&114858&3790336&3353460&91848366\\
159&3240&139590&117018&3885301&3437531&94727319\\
160&3281&142209&119205&3982229&3523243&97677964\\
161&3321&144837&121419&4080903&3610461&100701803\\
162&3362&147538&123661&4181362&3699361&103800552\\
163&3403&150248&125929&4283578&3789835&106975060\\
164&3445&152967&128226&4387645&3882032&110226672\\
165&3486&155742&130550&4493778&3975811&113556854\\
166&3528&158544&132902&4601761&4071354&116966891\\
167&3570&161394&135282&4711691&4168549&120458340\\
\hline\caption{Values of $b_{max}(z_l,s_{max})$ and
$b_{max}^*(z_l,s_{max})$}\label{t2}
\end{longtable}}
\end{center}

\end{document}